\documentclass{elsarticle}

\biboptions{sort&compress,numbers}
\usepackage[nodots]{numcompress}
\usepackage{lineno,url}
\usepackage{epstopdf}
\usepackage{amsmath,amsfonts,amssymb}
\usepackage{amssymb,amsmath,amsthm}
\usepackage{overpic}
\usepackage{contour}\contourlength{0.7pt}
\usepackage{tikz}\usetikzlibrary{shapes,arrows}
\usepackage{versions}
\usepackage{natbib}
\usepackage{algorithm, algorithmic}
\usepackage{multirow}

\newtheoremstyle{mytheorem}{5pt plus 5pt minus 3pt}{4pt plus 3pt minus 1.5pt}
	{\itshape}{}{\bfseries}{.}{1ex plus 1ex minus .5ex}{}
\newtheoremstyle{mydef}{5pt plus 5pt minus 3pt}{4pt plus 3pt minus 1.5pt}
	{}{0pt}{\bfseries}{.}{1ex plus 1ex minus .5ex}{}
\newtheoremstyle{myremark}{5pt plus 5pt minus 3pt}{4pt plus 3pt minus 1.5pt}
	{}{0pt}{\itshape}{.}{1ex plus 1ex minus .5ex}{}
\theoremstyle{mytheorem}
\newtheorem{prop}{Proposition}[section]

\theoremstyle{mydef}

\newtheorem{ex}[prop]{Example}
\theoremstyle{myremark}
\newtheorem{rem}{Remark}

\newcommand{\bb}{\mathbf b}

\newcommand{\de}{\mathrm{d}}

\newcommand{\bm}{\mathbf m}
\newcommand{\bx}{\mathbf x}
\newcommand{\D}{\mathbf D}
\newcommand{\DDt}{\widetilde{\mathbf D}}
\newcommand{\Q}{\mathcal Q}
\newcommand{\Qt}{\widetilde{\Q}}
\newcommand{\taut}{\widetilde{\tau}}
\newcommand{\omegat}{\widetilde{\omega}}
\newcommand{\tf}{\footnotesize}
\newcommand{\br}{\mathbf r}
\def\R {\mathbb R}
\newcommand{\FF}{\mathcal F}
\newcommand{\FFt}{\widetilde{\mathcal F}}

\newcommand{\Om}{\Omega}
\newcommand{\Omt}{\widetilde{\Omega}}

\newcommand{\XX}{\mathcal X}
\newcommand{\XXt}{\widetilde{\mathcal X}}
\newcommand{\xxt}{\widetilde{\bx}}
\newcommand{\bmt}{\widetilde{\bm}}
\newcommand{\mt}{\widetilde{m}}
\newcommand{\xt}{\widetilde{x}}
\newcommand{\St}{\widetilde{S}}
\newcommand{\Dt}{\widetilde{D}}
\newcommand{\card}{\textnormal{card}}

\def\Acaption#1#2{\caption{#2}\vspace*{-#1}}

\definecolor{Mgreen}{RGB}{34,139,34}
\definecolor{blau}{rgb}{0.15,0.2,0.5}
\definecolor{gray}{rgb}{0.5,0.5,0.5}
\definecolor{drot}{rgb}{0.7,0,0.1}
\definecolor{gelb}{rgb}{.55,.40,.1}
\definecolor{magenta}{rgb}{1.,0.,1.}
\definecolor{cyan}{rgb}{0.,1.,1.}
\definecolor{green}{rgb}{0.,1.,0.}
\definecolor{Morange}{rgb}{1.,0.5,0.}

\usepackage{array}
\newcolumntype{C}[1]{>{\centering\let\newline\\\arraybackslash\hspace{0pt}}m{#1}}

\begin{document}

\begin{frontmatter}

\title{Optimal quadrature rules for isogeometric analysis}

\author[NumPor]{Michael Barto\v{n}\corref{cor1}}
\ead{Michael.Barton@kaust.edu.sa}
\author[NumPor,AMCS]{Victor Manuel Calo}
\ead{Victor.Calo@kaust.edu.sa}

\cortext[cor1]{Corresponding author}

\address[NumPor]{Center for Numerical Porous Media, King Abdullah University of Science and Technology, Thuwal 23955-6900, KSA}%
\address[AMCS]{Applied Mathematics $\&$ Computational Science and Earth Science $\&$ Engineering, \\
King Abdullah University of Science and Technology, Thuwal 23955-6900, KSA}

\begin{abstract}
We introduce optimal quadrature rules for spline spaces that are frequently used
in Galerkin discretizations to build mass and stiffness matrices.
Using the homotopy continuation concept \cite{Homotopy-2015} that transforms
optimal quadrature rules from source spaces to target spaces, we derive optimal rules for splines
generated above finite domains.
Starting with the classical Gaussian quadrature for polynomials, which is an optimal rule for a discontinuous
odd-degree space, we derive rules for target spaces of higher continuity.
We further show how the homotopy methodology handles cases where the source and target rules
require different numbers of optimal quadrature points.
We demonstrate it by deriving optimal
rules for various odd-degree spline spaces, particularly with non-uniform
knot sequences and non-uniform multiplicities.
We also discuss convergence of our rules to their asymptotic counterparts, that is,
the analogues of the midpoint rule of Hughes et al. \cite{Hughes-2010}, that are exact
and optimal for infinite domains.
For spaces of low continuities, we numerically show that the derived rules quickly converge to their asymptotic
counterparts as the weights and nodes of a few boundary elements differ from the asymptotic values.
\end{abstract}


\end{frontmatter}

\section{Introduction and motivation}\label{intro}


Numerical integration is a basic ingredient of Galerkin discretizations,
which are at the core of finite elements and isogeometric analyses.
We aim to make this fundamental building block of computation more economical
by introducing new optimal quadrature rules for numerical integration.

We derive optimal rules that result in important computational savings and can
be used in many applied fields. Due to the popularity and success of isogeometric analysis
in many relevant engineering applications \cite{OptimizedWings-2014,DynEarthSol2D-2013,Bazilevs-2010,collier2014computational,cortes2015performance,vignal2015coupling,
auricchio2010isogeometric,wall2008isogeometric,niemi2012isogeometric,borden2011isogeometric,duddu2012finite,motlagh2013simulation,
calo2008multiphysics,calo2008simulation,gao2014fast,chang2012isogeometric,bazilevs2007weak,gomez2010isogeometric,
gomez2008isogeometric,bazilevs2008isogeometric,vignal2013phase,lipton2010robustness,temizer2011contact,hsu2011high,elguedj2008projection},
we focus on spline spaces, which are piece-wise polynomial functions with controlled global smoothness \cite{Piegl-1976,IGA-2009}.
Thus, we derive optimal (Gaussian) quadrature rules for spline spaces of arbitrary order and continuity.


A spline space is uniquely determined by its \emph{knot vector}. This knot vector is
a non-decreasing sequence of real numbers called \emph{knots} and the multiplicity of each knot determines the smoothness
of the polynomial pieces at the knot location.
A detailed introduction on splines can be found, e.g., in \cite{deBoor-1972,Hoschek-2002-CAGD,Elber-2001}.


We call a \emph{quadrature rule}, or simply a \emph{quadrature}, an \textit{$m$-point rule},
if $m$ evaluations of a function $f$ are needed to approximate its weighted integral
over a closed interval $[a,b]$
\begin{equation}\label{eq:GaussQuad}
\int_a^b w(x) f(x) \, \mathrm{d}x = \sum_{i=1}^{m} \omega_i f(\tau_i) + R_{m}(f),
\end{equation}
where $w$ is a fixed non-negative \emph{weight function} defined over $[a,b]$.
The rule is required to be \emph{exact}, i.e., $R_m(f) \equiv 0$
for each element of a spline space $S$.
A rule is \emph{optimal} if $m$ is the minimum number of
\emph{weights} $\omega_i$ and \textit{nodes} $\tau_i$ (points at which $f$ is evaluated).

For discontinuous spline spaces (polynomials), the optimal rule is known to be the classical Gaussian quadrature \cite{Gautschi-1997}
with the \emph{order of exactness} $2m-1$, that is, only $m$ evaluations are needed to exactly integrate
any polynomial of degree at most $2m-1$. Consider a sequence of polynomials $(q_0,q_1,\ldots,q_m,\ldots )$  that form an orthogonal basis
with respect to the scalar product
\begin{equation}
<f,g> = \int_a^b f(x)g(x)w(x) \de x.
\end{equation}
The quadrature points are the roots of the $m$-th orthogonal polynomial $q_m$ which
in the case when $w(x)\equiv 1$ is the degree-$m$ Legendre polynomial \cite{Szego-1936}.



The quadrature rules for splines have been studied since the late 50's \cite{Schoenberg-1958,Micchelli-1972,Micchelli-1977}.
Micchelli and Pinkus \cite{Micchelli-1977} proved that, for spaces with uniform continuities, i.e. knots' multiplicities,
there always exists an optimal quadrature formula with the following number of necessary evaluations:
\begin{equation}\label{eq:Micchelli}
     d + 1 + i = 2 m,
\end{equation}
where $d$ is the polynomial degree, $i$ is the total number of interior knots (when counting multiplicities),
and $m$ is the number of optimal nodes.

Regarding concrete optimal quadrature rules over a finite domain, very little is known in the literature.
Even though \cite{Micchelli-1977} gives a range of knot spans,
among which each particular node lies, the set of subintervals is too large to even build the corresponding polynomial system.
Therefore, formulating it as an optimization problem is rather difficult as it is highly non-linear and a good initial guess in essential. 
To the best of our knowledge, however, there is no theory that tightly bounds each node to give a good initialization for the numerical optimization.

For $C^1$ cubic splines with uniform knots, Nikolov \cite{Nikolov-1996} derives a recursive formula that starts at the boundary of the interval
and parses towards its middle, recursively computing all the optimal nodes and weights. Computationally, this result is favourable
as it computes the nodes and weights in a \emph{closed form}, i.e.,
without the need of any numerical solver \cite{SolverGershon-2001,Patrikalasis-1993,SolverET-2012}.
We have generalized this result to a special \emph{class of spaces} built above non-uniform knot vectors \cite{Quadrature31-2014},
called symmetrically stretched B-splines. 
Recently, we have described a similar recursive pattern as Nikolov for $C^1$ quintic spline spaces with uniform knots \cite{Quadrature51-2014}
and derived optimal rules for them.
Interestingly, the optimal rules are still explicit, even though five degree polynomials are involved.
We proved that there exists an algebraic factorization in every
step of the recursion which makes the rule explicit and therefore computationally cheap.

Recently, we demonstrated a relation between Gaussian quadrature rules for $C^1$ and $C^2$ cubic splines \cite{Homotopy-2015}.
We showed that for particular pairs of spaces that require the same number of optimal nodes, the quadrature rule can be ``transferred''
from the \emph{source} space to the \emph{target} space, preserving the number of internal knots in (\ref{eq:Micchelli})
and therefore the number of optimal quadrature points.
This transfer relies on the application of homotopy continuation to preserve the algebraic structure of the system
as the spline space evolves \cite{Wampler-2005}.
Starting with a known optimal rule for a $C^1$ cubic spline space, the source knot vector $\XXt$ (with double knots)
gets transformed to the desired configuration $\XX$ (single knots). The quadrature rule is considered as a high-dimensional point,
and satisfies a certain well-constrained system of polynomial equations. These equations reflect the fact that the rule must exactly integrate
all the basis functions defined by $\XXt$. Using homotopy continuation, the rule is numerically traced as $\XXt$,
and consequently the above-built spline space, continuously transforms to its target configuration $\XX$.

For higher degrees than five, we expect to require a numerical solver. But more importantly,
to prove the correct layout of nodes (exact pair of knots that tightly delimit every optimal node) does not seem to be feasible at all.
From this perspective, homotopy continuation based
methods \cite{Homotopy-2015} are a feasible approach because the correct layout of nodes comes automatically
from the application of the continuation methodology.

In this paper, we further develop the homotopy continuation methodology and derive optimal rules for odd-degree spline spaces.
Whereas in the case of $C^2$ cubic splines defined above an odd number of elements
there always exists a $C^1$ cubic spline space that possesses the same number of optimal quadrature nodes,
this is not true in general. For example, for $C^{-1}$ and $C^1$ septic spline spaces the counting of degrees of freedom
is non-trivial, while optimal rules exist for each particular space.
The number of necessary evaluations that satisfy (\ref{eq:Micchelli}) is
a strict constraint on the number of source and target elements.
However, we show that even for such spaces, the transformation of Gaussian rules is possible
for properly chosen element counts.
That is, we show that the homotopy concept can be generalized to handle scenarios,
where the source and target spaces require different numbers of optimal nodes to exactly integrate each of them.

The rest of the paper is organized as follows. Section~\ref{sec:matrices} describes the particular class of spline spaces
for which the optimal rules are derived, Section~\ref{sec:splines} summarizes a few basic properties of spline spaces of odd degree
and shows possible continuous transformations between them. In Section~\ref{sec:homotopy}, we show how
the homotopy continuation can be used to derive new optimal quadrature rules.
Section~\ref{sec:ex} discusses the results and the validity of the new quadratures obtained.
We conclude the paper discussing the results, their relevance and future research directions.


\section{Mass and stiffness matrices for isogeometric analysis}\label{sec:matrices}

We provide Gaussian quadrature rules for a class of spline spaces
that arise when solving elliptic partial differential equations (PDEs).
The key idea is to find the smallest spline space that contains all
terms appearing when building the mass and stiffness matrices.
For each of these spaces, we aim to derive the appropriate Gaussian quadrature rule
that accounts for arbitrary knot spacings as well as multiplicities.

Consider a spline space of degree $p$ and continuity $k$, $S_{p,k}$, then the Grammian (mass) matrix contains
elements from $S_{2p,k}$. When differentiating once $S_{p,k}$, one obtains $S_{p-1,k-1}$ and, in general
when applying an $l$-th order differential operator, the spline space is $S_{p-l,k-l}$ and the $L_2$
scalar product of its elements lies in $S_{2(p-l),k-l}$.
With $l$ differentiations, we obtain a hierarchy of spline spaces as shown in Fig.~\ref{fig:hierarchy}
and the minimal space that contains terms of both mass and stiffness matrices is $S_{2p,k-l}$.
The first limitation stems in the fact that this space is of even degree.
Nevertheless, as is commonly done in finite element theory, we use Gauss rules
that are optimal for the smallest higher odd-degree polynomials \cite{Ciarlet-2002,Brenner-2009,Strang-2008}.
Therefore, we focus on rules for spline spaces of odd degree $d:=2p+1$ and continuity $c:=k-l$ for which we derive optimal rules.

 \begin{figure}[!tb]
\vrule width0pt\hfill
\begin{equation*}
\renewcommand{\arraystretch}{1.4}
\begin{array}{ccc}
(p,k)     & \rightarrow & (2p,k)   \\
(p-1,k-1) & \rightarrow & (2(p-1),k-1)   \\
\vdots    & & \vdots \\
(p-l,k-l) & \rightarrow & (2(p-l),k-l)
\end{array}
\end{equation*}
 \hfill \vrule width0pt\\
 \vspace{-25pt}
\Acaption{1ex}{A hierarchy of spline spaces used when building mass and stiffness matrices.
Left: the original spline space of degree $p$ and continuity $k$ is being differentiated $l$-times.
Right: the corresponding spline spaces that contain the appropriate scalar products for Galerkin methods 
are shown.}\label{fig:hierarchy}
 \end{figure}


%
%


\begin{ex}
Let $p=3$, $k=2$, $l=1$, then we have
\begin{equation*}
\renewcommand{\arraystretch}{1.4}
\begin{array}{ccc}
(6,2)     & \subset & (6,1)   \\
\cup      &  & \cup   \\
(4,2)     & \subset & (4,1)
\end{array}
\end{equation*}
The inclusion relations follow directly from the fact that the corresponding knot vectors
are nested \cite{deBoor-1978}. In the context of finite elements in 1D,
the elements of the mass matrix belong to the $(6,2)$-space
whilst the scalar products that fill the stiffness matrix belong to $(4,1)$-space.
To minimize the computational cost and simplify implementation, we choose an optimal
rule that can integrate all terms. Thus, we seek the smallest space that contains both
$(4,1)$ and $(6,2)$, and is of odd degree, i.e, $(d,c)=(7,1)$.
\end{ex}


\section{Continuous transformations of spline spaces}\label{sec:splines}

We recall several properties of spline spaces. Consider a knot vector
\begin{equation}\label{eq:x}
\begin{array}{ccccccc}
\XX_N =  & (a= &\underbrace{x_0, \dots,x_0,} & \underbrace{x_1, \dots, x_{1},} & \dots & \underbrace{x_N,\dots,x_N}& =b)\\
& & m_0 & m_1 & & m_N &
\end{array}
\end{equation}
and for the sake of simplicity, we split $\XX_N:=(\bx,\bm)$ into the \emph{domain partition} $\bx$, $\bx \in \mathbb R^N$,
and the \emph{vector of multiplicities} $\bm$, $\bm \in \mathbb N^N$, and write
\begin{equation}\label{eq:split}
\bx =  (x_0, \dots,x_N), \quad \quad \bm = (m_0,\dots, m_N).
\end{equation}
We further recall $1\leq m_i \leq d+1$, $i=0,\dots,N$ and assume $\XX_N$ is an \emph{open knot} vector on $[a,b]$, that is, $m_0=m_N = d+1$,
and $d$ is odd.
We denote by $\pi_d$ a space of polynomials of degree at most $d$ and define
the spline space associated to $\XX_N$ as
\begin{equation}\label{eq:targetSpace}
S_{\bx,\bm}^{N,d} = \{ f\in C^{d-m_k} \,\, \textnormal{at} \, x_k, k=0,\dots,N \,\, \textnormal{and} \,  f|_{(x_{k-1},x_{k})} \in \pi_d, k=1,\dots,N\}.
\end{equation}
Our goal is to derive a Gaussian rule for this target space $S_{\bx,\bm}^{N,d}$.
To do so, we define an associated source space for which the optimal rule is known. In particular we consider a
uniform knot vector with maximum knot multiplicities
\begin{equation}\label{eq:xt}
\begin{array}{ccccccc}
\XXt_n =  & (a= &\underbrace{\xt_0, \dots,\xt_0,} & \underbrace{\xt_1, \dots, \xt_{1}},& \dots & \underbrace{\xt_n,\dots,\xt_n}& =b)\\
& & d+1 & d+1 & & d+1 &
\end{array}
\end{equation}
and using analogous notation to (\ref{eq:split}), i.e., $\XXt_n:=(\xxt,\bmt)$, we obtain
\begin{equation}\label{eq:splitt}
\xxt =  (\xt_0, \dots, \xt_n), \quad \bmt = (\mt_0,\dots, \mt_n),
\end{equation}
where $h=\frac{b-a}{n}$, $\xt_i = a + hi$ and $\mt_i = d+1$ for $i=1,\dots,n$.
We define the \emph{source} spline space as
\begin{equation}\label{eq:sourceSpace}
\St_{\xxt,\bmt}^{n,d} = \{ f\in C^{-1} \,\, \textnormal{at} \, \xt_k, k=0,\dots,n \,\, \textnormal{and} \,  f|_{(\xt_{k-1},\xt_{k})} \in \pi_d, k=1,\dots,n\}.
\end{equation}
Our source space $\St_{\xxt,\bmt}^{n,d}$ is chosen to be discontinuous, that is,
$C^{-1}$ continuous space of polynomials because for this space, when $d$ is odd,
the classical Gaussian quadrature for polynomials can be used on every element to derive an optimal source rule.
An immediate question that arises is how
to set $n$ such that the number of nodes is optimal for our target space (\ref{eq:targetSpace}).
To answer this question, we further define the cardinality of the knot vectors as
\begin{equation}\label{eq:Card}
\card(\XX_N) = \sum_{j=0}^N m_i, \quad \textnormal{and} \quad \card(\XXt_n) =\sum_{i=0}^n \mt_i
\end{equation}
and require
\begin{equation}\label{eq:TotalKnots}
\card(\XX_N) \leq \card(\XXt_n). 
\end{equation}
We want to transform $\St_{\xxt,\bmt}^{n,d}$ into $S_{\bx,\bm}^{N,d}$, so ideally, we wish to have the same number of knots.
In such a case, the spaces have the same dimension and one can apply the scheme described in \cite{Homotopy-2015} where the
vectors of the same cardinality were transformed into one another, yielding an optimal rule.

Unfortunately, finding an $n$ such that the equality in (\ref{eq:TotalKnots}) holds is not always possible
for an arbitrary choice of $\card(\XX_N)$.
Nevertheless, one can always build a source space $\St_{\xxt,\bmt}^{n,d}$ such that the total number of
knots, $\card(\XXt_n)$, exceeds $\card(\XX_N)$ by at most $r \leq d-1$ knots.
That is, $\St_{\xxt,\bmt}^{n,d}$ and $S_{\bx,\bm}^{N,d}$ have different dimensions,
with the source having the minimum number of extra basis functions.

\begin{figure}[!tb]
\vrule width0pt\hfill
 \begin{overpic}[width=0.49\textwidth,angle=0]{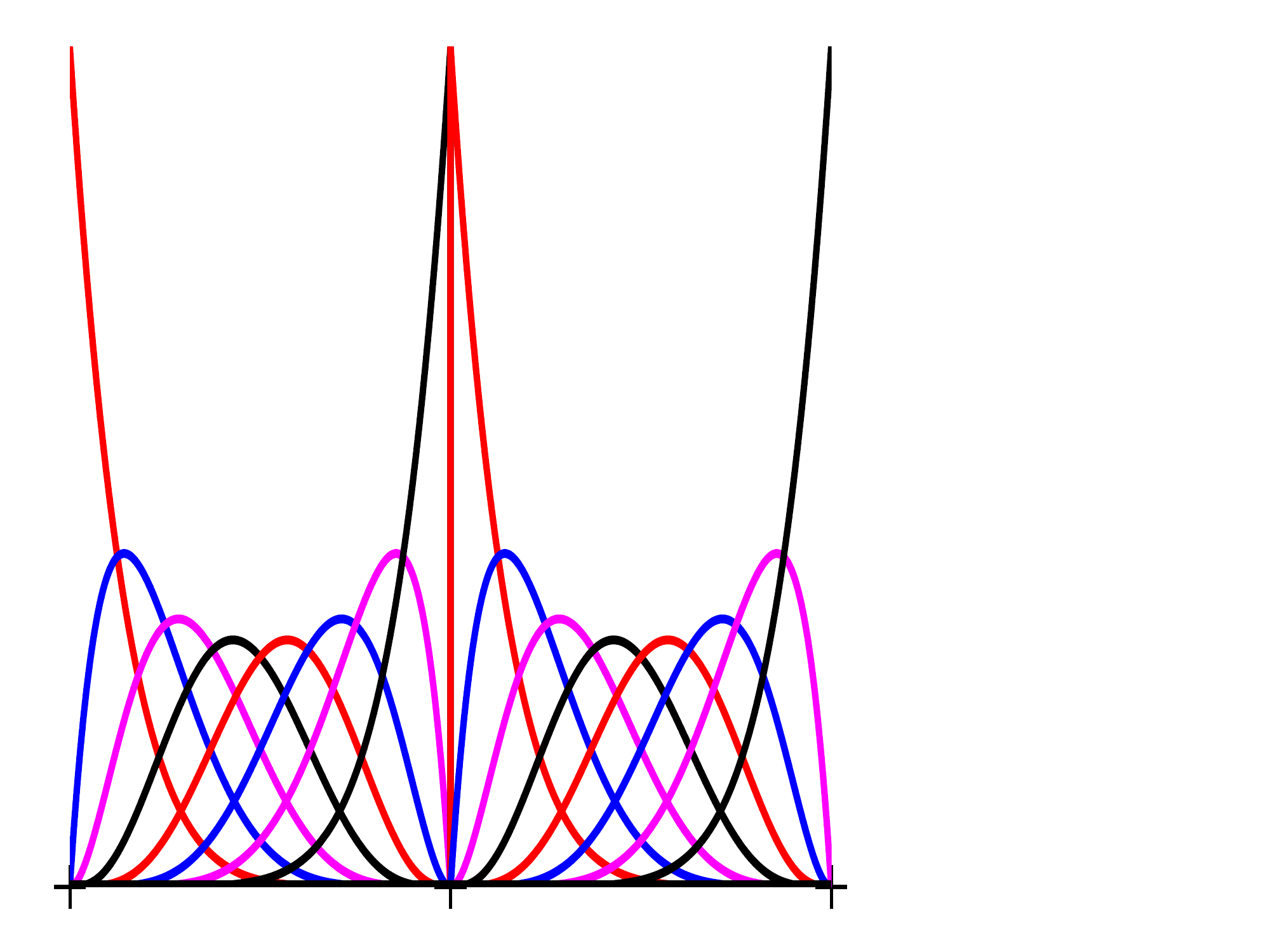}
 \put(10,60){\fcolorbox{gray}{white}{$\St_{\xxt,\bmt}^{2,7}$}}
 \put(38,58){\scriptsize$\xxt = (0,\frac{1}{2},1)$}
 \put(38,64){\scriptsize$\bmt = (8,8,8)$}
 \put(69,64){\scriptsize$\dim(\St_{\xxt,\bmt}^{2,7}) = 16$}
    \put(5,-2){\scriptsize$8$}
    \put(2,7){\scriptsize$a$}
    \put(30,-2){\scriptsize$8$}
    \put(65,-2){\scriptsize$8$}
    \put(70,7){\scriptsize$b$}
	\end{overpic}
 \hfill
 \begin{overpic}[width=0.49\columnwidth,angle=0]{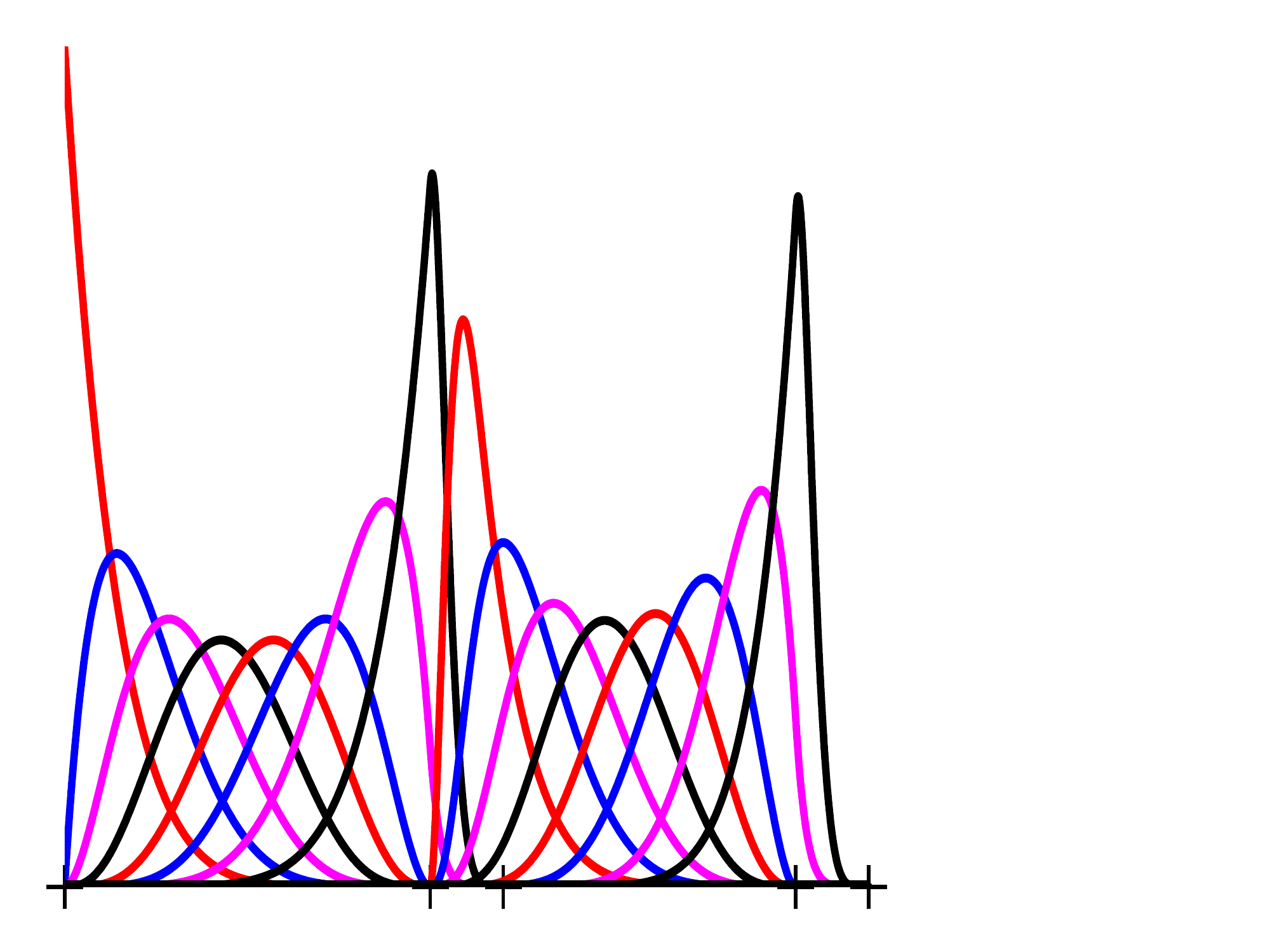}
 \put(-15,30){\huge$\rightarrow$}
  \put(35,63){\scriptsize$\bx = (0,\frac{1}{2},\frac{3}{5},1,\frac{11}{10})$}
  \put(35,69){\scriptsize$\bm = (8,6,2,6,2)$}
  \put(5,-2){\scriptsize$8$}
  \put(33,-2){\scriptsize$6$}
  \put(39,-2){\scriptsize$2$}
  \put(62,-2){\scriptsize$6$}
  \put(67,-2){\scriptsize$2$}
	\end{overpic}\hfill \vrule width0pt\\[1ex]
\vrule width0pt\hfill
 \begin{overpic}[width=.49\columnwidth,angle=0]{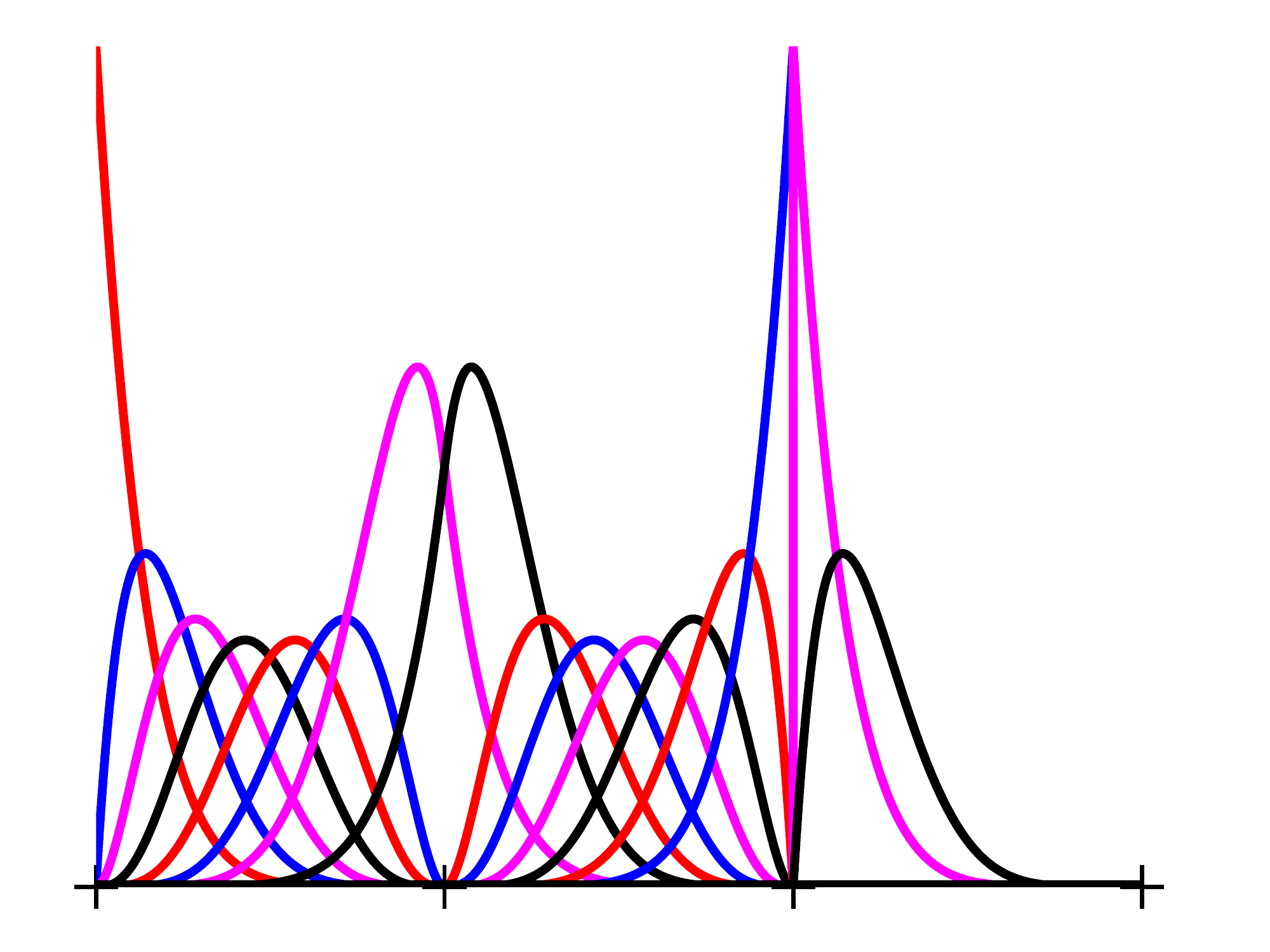}
 \put(10,65){\fcolorbox{gray}{white}{$S_{\bx,\bm}^{2,7}$}}
 \put(30,60){\scriptsize$\bx = (0,\frac{1}{2},1,\frac{3}{2})$}
 \put(30,66){\scriptsize$\bm = (8,6,8,2)$}
 \put(66,66){\scriptsize$\dim(S_{\bx,\bm}^{2,7}) = 14$}
    \put(5,-2){\scriptsize$8$}
  \put(32,-2){\scriptsize$6$}
  \put(3,7){\scriptsize$a$}
  \put(65,7){\scriptsize$b$}
  \put(60,-2){\scriptsize$8$}
  \put(89,-2){\scriptsize$2$}
  \put(66,40){\scriptsize$D_{15}$}
  \put(74,16){\scriptsize$D_{16}$}
	\end{overpic}\hfill
 \begin{overpic}[width=.49\textwidth,angle=0]{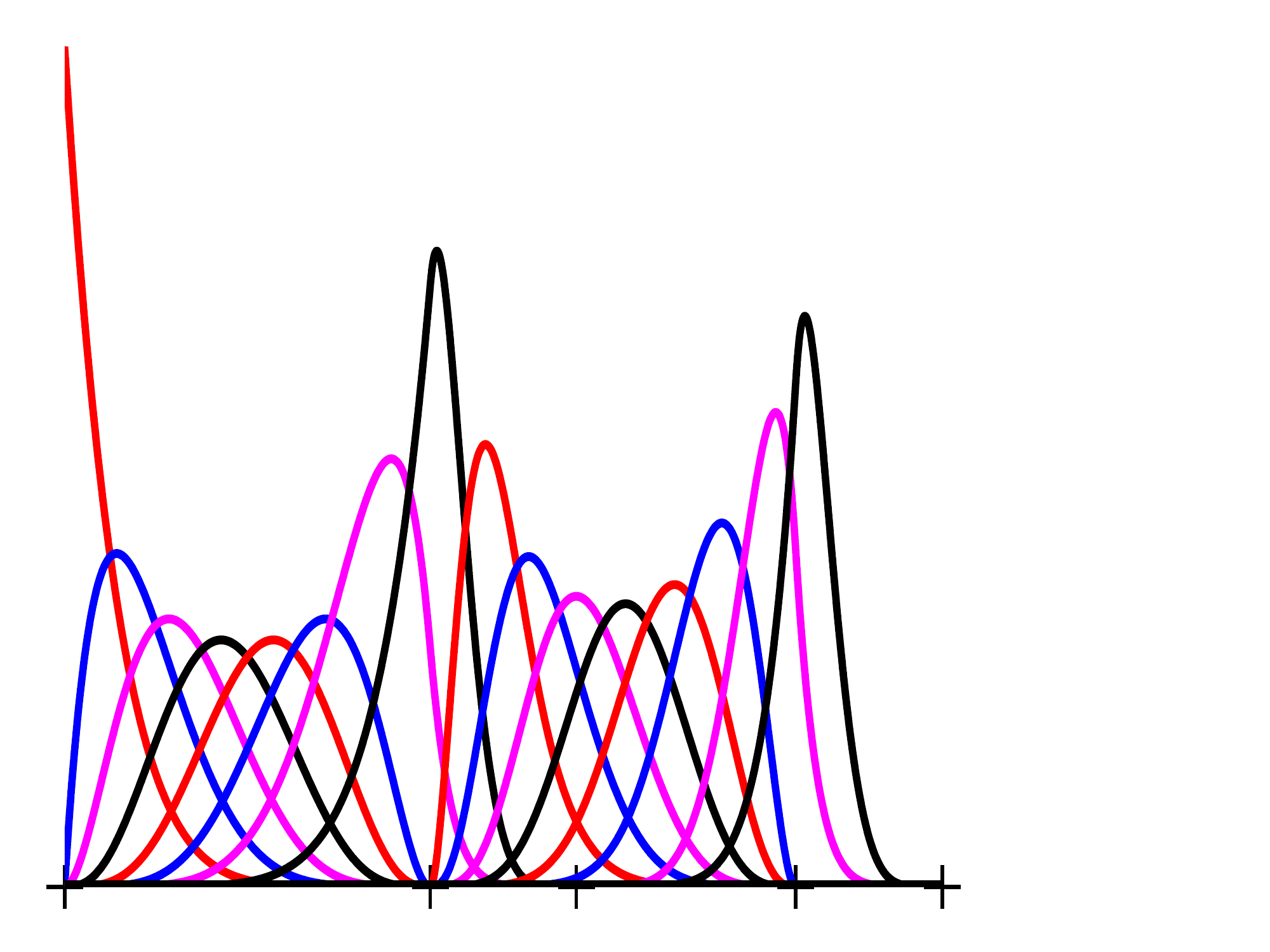}
 \put(47,60){\huge$\downarrow$}
  \put(5,-2){\scriptsize$8$}
  \put(33,-2){\scriptsize$6$}
  \put(45,-2){\scriptsize$2$}
  \put(62,-2){\scriptsize$6$}
  \put(75,-2){\scriptsize$2$}
	\end{overpic}\hfill \vrule width0pt\\[1ex]
\vrule width0pt\hfill
 \begin{overpic}[width=.49\columnwidth,angle=0]{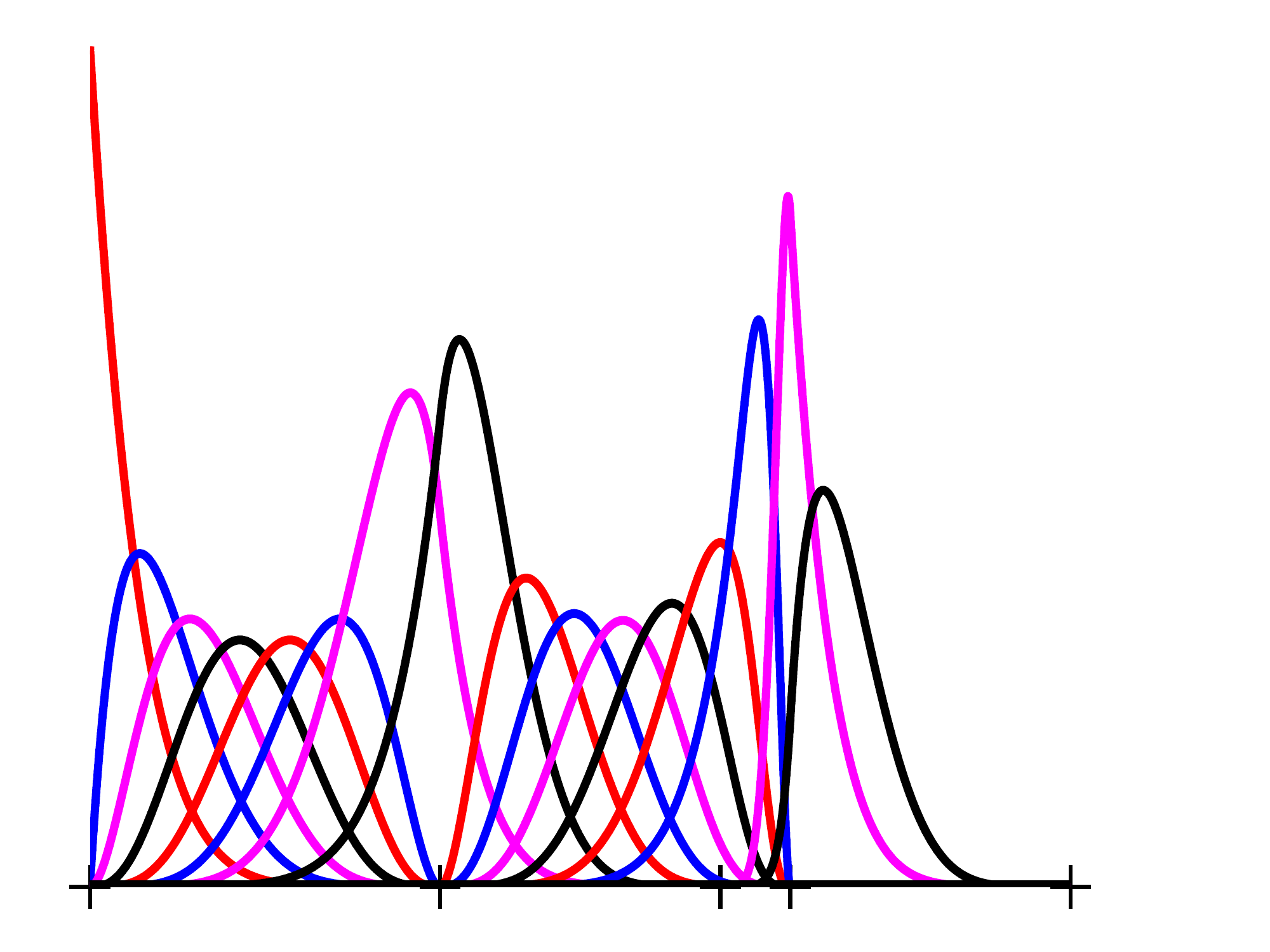}
\put(47,60){\huge$\uparrow$}
  \put(5,-2){\scriptsize$8$}
  \put(34,-2){\scriptsize$6$}
  \put(55,-2){\scriptsize$2$}
  \put(62,-2){\scriptsize$6$}
  \put(83,-2){\scriptsize$2$}
	\end{overpic}\hfill
 \begin{overpic}[width=.49\textwidth,angle=0]{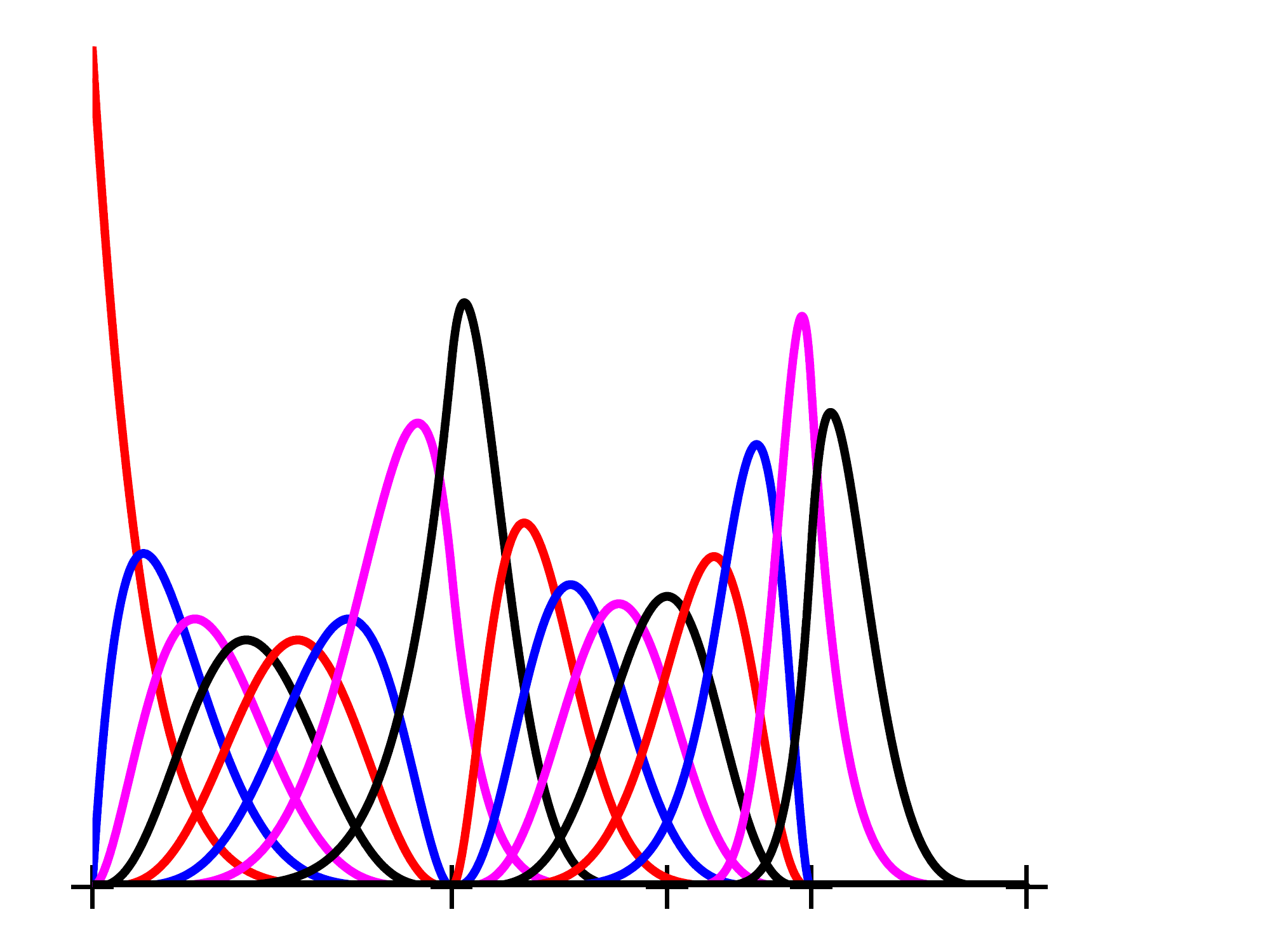}
\put(47,60){\huge$\downarrow$}
 \put(-10,30){\huge$\leftarrow$}
  \put(5,-2){\scriptsize$8$}
  \put(35,-2){\scriptsize$6$}
  \put(51,-2){\scriptsize$2$}
  \put(63,-2){\scriptsize$6$}
  \put(81,-2){\scriptsize$2$}
	\end{overpic}\hfill \vrule width0pt\\[-2ex]
\Acaption{1ex}{Continuous transformation of a spline spaces of degree $d=7$ on an interval $[a,b]$.
The source space over two elements ($n=2$) with $C^{-1}$ continuity, $\St^{2,7}_{\xxt,\bmt}$,
is being transformed to the $C^{1}$-continuous space, $S^{2,7}_{\bx,\bm}$, with $N=2$ elements.
  The transformation is realized via a continuous knot vector modification $\XXt_n \rightarrow \XX_N$.
  The number bellow each knot in the figure refers to its multiplicity.
   During the continuation, two boundary knots are pulled outside the domain $[a,b]$
   and simultaneously two internal knots are pulled towards $b$, which results in $D_{15}$ and $D_{16}$
   loosing their support over $[a,b]$ in the limit (middle left). This reduces the dimension of the spline space
   that acts on $[a,b]$ from $16$ to $14$.
   The source space requires $8$ Gaussian nodes while the target space needs only $7$ nodes on $[a,b]$, which is in accordance
   with the optimal count of Micchelli and Pinkus (\ref{eq:Micchelli}).}\label{fig:Trans7}
 \end{figure}

We seek optimal rules for target spaces of uniform continuity $c$, we have
\begin{equation}\label{eq:C1source}
\bm = (d+1,d-c,\dots, d-c,d+1)
\end{equation}
and a direct computation of the dimensions of the source and target space gives
%
\begin{equation}\label{eq:dims}
\renewcommand{\arraystretch}{1.3}
\begin{array}{lcl}
\dim(\St_{\xxt,\bmt}^{n,d}) & = & n(d+1),  \\
\dim(S_{\bx,\bm}^{N,d}) & = & Nd - Nc + c + 1.
\end{array}
\end{equation}
With the sharp inequality in (\ref{eq:TotalKnots}), the Gaussian quadrature in the source space
requires more quadrature nodes than the optimal rule for the target space.
The key idea is, however, that the $r$ extra knots can be moved outside $[a,b]$ which will force $r$ basis functions
to loose their support on $[a,b]$, see Fig.~\ref{fig:Trans7}.

We further define $\DDt = \{\Dt_i\}_{i=1}^{n(d+1)}$, the basis of $\St^{n,d}_{\xxt,\bmt}$, as
\begin{equation}\label{eq:Bernstein}
\renewcommand{\arraystretch}{1.2}
\begin{array}{lcl}
\Dt_{(d+1)k-d}(t) & = & [\xt_{k-1},\xt_{k-1},\dots, \xt_{k-1},\xt_k](. - t)_{+}^{d}   \\
 &\vdots & \\
\Dt_{(d+1)k-1}(t) & = & [\xt_{k-1},\xt_{k-1},\xt_{k},\dots,\xt_k](. - t)_{+}^{d}   \\
\Dt_{(d+1)k}(t)   & = & [\xt_{k-1},\xt_{k},\dots, \xt_{k},\xt_k](. - t)_{+}^{d}, \\
\end{array}
\end{equation}
where $[.]f$ stands for the divided difference and $u_{+} = \max(u,0)$ is the truncated
power function and $k=1,\dots, n$.
Because of the maximum multiplicity $d+1$ at each knot of the knot vector $\XXt$, see (\ref{eq:xt}), $\Dt_k$ are the Bernstein basis functions
of degree $d$, see Fig.~\ref{fig:Trans7} top left, and therefore
\begin{equation}\label{interiorIntegral}
I[\Dt_i] = \frac{1}{d+1}\; \textnormal{for} \quad i = 1,2,\ldots,n(d+1),
\end{equation}
where $I[f]$ stands for the integral of $f$ over the interval $[a,b]$, see e.g. \cite{Hoschek-2002-CAGD}.


For the $C^c$-continuous target source space $S_{\bx,\bm}^{N,d}$, i.e., for setting $\bm$ as in (\ref{eq:C1source}),
the basis $\D = \{D_i\}_{i=1}^{Nd - Nc + c + 1}$ is built analogously to (\ref{eq:Bernstein}) with the difference that the internal
knots possess multiplicity $d-c$ instead of $d+1$. However, as already mentioned,
the dimensions in general do not match
\begin{equation}\label{eq:dimmismatch}
Nd-Nc+c+1 \leq n(d+1),
\end{equation}
but we still want to find a continuous transformation of $\St_{\xxt,\bmt}^{n,d}$ into $S_{\bx,\bm}^{N,d}$
over $[a,b]$ which is governed by
\begin{equation}\label{eq:Trans}
\XXt_n \rightarrow \XX_N.
\end{equation}
Our aim is to transfer the optimal rule from $\St_{\xxt,\bmt}^{n,d}$ to $S_{\bx,\bm}^{N,d}$.

We show now how to set the source space such that
the transformation to the target space
admits the optimal quadrature node count (\ref{eq:Micchelli}) on $[a,b]$.
For example, for $c=1$,
to satisfy the desired inequality of the dimensions in (\ref{eq:dimmismatch}), we set
\begin{equation}\label{eq:n}
n := \left \lceil{\frac{Nd-N+2}{d+1}}\right \rceil
\end{equation}
where $\left \lceil{\cdot}\right \rceil$ is the ceiling function
and therefore $n$ is the smallest integer that satisfies (\ref{eq:dimmismatch}).
The dimension difference $r = \dim(\St_{\xxt,\bmt}^{n,d}) - \dim(S_{\bx,\bm}^{N,d})$ is even
because both sides of (\ref{eq:dimmismatch}) are even for odd $d$. Consequently
\begin{equation}\label{eq:r}
r = 0, 2, \dots, d-1.
\end{equation}



\begin{rem}
For a general vector of multiplicities of the target space, one has to adapt
$n$ in (\ref{eq:n}) accordingly.
The parity of $r$ stays even for any target spaces of constant odd continuity $C^c$,
because the right-hand side of (\ref{eq:dimmismatch}) is
even for any odd $c$. For even $c$, one needs additionally $N$ to be odd, in order to be in accordance with
the optimal count of Micchelli and Pinkus \cite{Micchelli-1977}.
In this work, we consider only spaces that meet optimal parity count in (\ref{eq:Micchelli}), i.e.,
require $m$ nodes to exactly integrate functions from a spline space of dimension $2m$.
%
\end{rem}

We have shown how to find the smallest $n$ such that the source space has larger or equal dimension than
the target one and that the difference of dimensions is always even. 
This parity allows us to derive optimal quadrature rules.
The difference $r$ is the number of \emph{redundant} basis functions
that do not to contribute to the target rule.
Therefore, the number of nodes that should become redundant during continuation is $\frac{r}{2}$.
As $r$ basis functions are forced to loose their support on $[a,b]$, by pulling $r$ boundary knots outside $[a,b]$
and $r$ internal knots to the same boundary, see Fig.~\ref{fig:Trans7}, $\frac{r}{2}$ quadrature points are forced to move to
the boundary whilst their weights vanish in the limit.
If $r$ is odd, one cannot expect to derive an optimal rule.

\begin{figure}[!tb]
\vrule width0pt\hfill
 \begin{overpic}[width=.5\textwidth,angle=0]{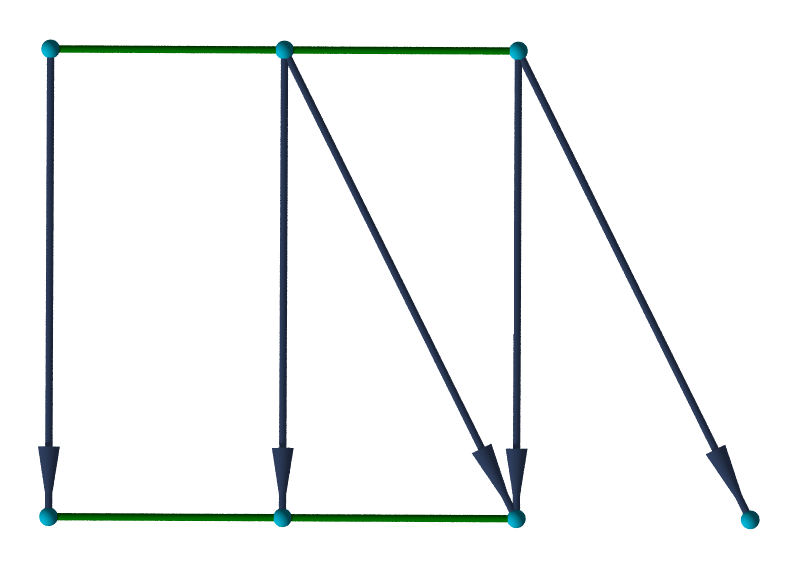}
    \put(-25,5){\fcolorbox{gray}{white}{$\XX_N$}}
    \put(-25,60){\fcolorbox{gray}{white}{$\XXt_n$}}
    \put(0,68){\small$a=\xt_{0}$}
    \put(60,68){\small$b=\xt_{2}$}
    \put(100,60){\small$\bmt=(8,8,8)$}
    \put(0,40){\small$8$}
    \put(30,40){\small$6$}
     \put(49,40){\small$2$}
     \put(59,40){\small$6$}
     \put(78,40){\small$2$}
     \put(100,5){\small$\bm=(8,6,8)$}
     \put(0,0){\small$a=x_{0}$}
    \put(60,0){\small$b=x_{2}$}
	\end{overpic}
 \hfill \vrule width0pt\\
 \vspace{-5pt}
\Acaption{1ex}{A target uniform knot vector $\XX_N$ for $N=2$ elements and an associated source knot vector $\XXt_n$, $n=2$, for
a septic spline space ($d=7$) defined over $[a,b]$ are shown.
The source knot vector exceeds the target one by $r=2$ knots, i.e., $\card(\XX_n) = 24$ whilst $\card(\XX_N) = 22$. The geodesic
knot transformation is visualized by the paths (dark blue), showing multiplicities of the transformed knots.}\label{fig:KnotsTrans}
 \end{figure}

We recall the fact that the external knots and the dimension differences between the source and target spaces
are both equal to $r$ and
the total number of knots of $\XX_N$ that lie inside $[a,b]$ is
\begin{equation}\label{eq:TotKnots}
K = Nd - Nc + c + r + 1.  
\end{equation}
Additionally, the right boundary $b$ has multiplicity $d+1$ in $\XX_N$.
Without loss of generality, we assume knot transformations where all $r$ boundary knots will leave the domain at $b$
and $r$ interior knots move from the inside to $b$, see Fig.~\ref{fig:KnotsTrans}.


There exist infinitely many knot transformations that satisfy (\ref{eq:Trans}). In particular, we consider transformations in which
every pair of mutually corresponding knots $(\XXt_n[i] \rightarrow \XX_N[i])$, $i=1,\dots, \card(\XX_N)$  follows the shortest path
between $\XXt_n$ and the corresponding subset of $\XX_N$ (with the same cardinality), when measured by the Euclidean metric in the knot
vector space. This path is linear in time and is called \emph{geodesic}, see \cite[page 5]{Homotopy-2015}.

There also exist infinitely many possible moves for the $r$ extra knots outside $b$.
To simplify the analysis, we assume that all the $r$
knots coincide at the target configuration, and we set the distance from $b$ as $\frac{1}{N}$.
That is, the boundary knot of multiplicity $r$ is also moved in a geodesic fashion
by the size of one element, see Fig.~\ref{fig:KnotsTrans}.

\section{Gaussian quadrature via homotopy continuation}\label{sec:homotopy}

In this section, we derive optimal quadrature rules for spline spaces with various polynomial degrees and continuities.
We derive optimal rules from other optimal ones as recently introduced in \cite{Homotopy-2015}
and refer the reader to that work for a more detailed description of the methodology.


\subsection{Gaussian quadrature}\label{ssec:guass}

We consider an odd-degree source space $\St_{\xxt,\bmt}^{n,d}$ with an optimal quadrature formula
\begin{equation}\label{quadratureS}
\Qt_a^b[f] = \sum_{i=1}^{m} \omegat_i f(\taut_i) = \int_{a}^{b} f(t) \mathrm{d}t, \quad f \in \St_{\xxt,\bmt}^{n,d}
\end{equation}
where $m = \frac{(d+1)n}{2}$, and the nodes are the $n$-tuples of the classical polynomial Gaussian points
computed on every element $[\xt_{i-1},\xt_{i}]$, $i=1,\dots, n$. This source rule is optimal and unique,
and assumes that each polynomial space is discontinuous at each element interface.

Consider the target space $S_{\bx,\bm}^{N,d}$ with the multiplicity vector defined in (\ref{eq:C1source}) and let $r$
be the even dimension difference
between the source and the target spaces, see (\ref{eq:r}). Then, according to \cite{Micchelli-1977},
there exists an optimal target rule
\begin{equation}\label{quadrature}
\Q_a^b[f] = \sum_{i=1}^{m-\frac{r}{2}} \omega_i f(\tau_i)= \int_{a}^{b} f(t) \mathrm{d}t, \quad f \in S_{\bx,\bm}^{N,d}.
\end{equation}

Thus, we have guaranteed existence of optimal rules for the source and target spaces, even though one space requires $m$ and
the latter only $m-\frac{r}{2}$ optimal nodes.
However, to the best of our knowledge, there are no theoretical results on existence of Gaussian quadrature for
spaces with non-uniform continuities. We therefore only conjecture that there also exists an optimal rule
for any intermediate knot vector and, consequently, for any intermediate spline space. Our numerical results, see Section~\ref{sec:ex},
support this claim.

As the source space gets transformed into the target space by continuously modifying the knot vector, via (\ref{eq:Trans}),
our goal is to trace the known source rule to derive an optimal rule for the target space.
Therefore $\Q$, represented by its nodes and weights, is a function of time $t$, $t\in [0,1]$.
To simplify notation, if no ambiguity is imminent,
we omit the time parameter and write $\tau_i$ instead of $\tau_i(t)$.
The source rule is $\Qt=\Q(0)$ and the target rule we wish to derive is $\Q=\Q(1)$.


\subsection{Homotopy continuation}\label{ssec:homotopy}

Polynomial homotopy continuation (PHC) is a numerical scheme commonly used to solve
polynomial systems of equations \cite{Wampler-1990,Wampler-2005}.
Given a polynomial system $\FF(\bx)=\mathbf 0$ that we want to solve,
the method uses the known roots of a simpler polynomial system (source) $\FFt(\bx)=\mathbf 0$
which is continuously transformed into the desired (target) solution.
We can therefore write
\begin{equation}\label{eq:SystemTime}
\FF(\bx,t)=\mathbf 0
\end{equation}
that at $t=0$ is the system whose roots we know, and at $t=1$ is the target system we aim to solve.
We refer the reader to \cite{Wampler-2005}
for a detailed explanation of PHC.
For the sake of completeness, we now briefly review the homotopic setting in the context of quadrature rules \cite{Homotopy-2015}:
the quadrature rule $\Q$ (the nodes and weights) is encoded as a high-dimensional point
$$
\bx = (\tau_1,\dots, \tau_m, \omega_1,\dots, \omega_m), \quad \bx \in \R^{2m},
$$
our source $2m \times 2m$ polynomial system $\FF(\bx,0) = \mathbf 0$ expresses that the source rule $\Qt$ (\ref{quadratureS})
exactly integrates the source basis $\DDt$, that is,
\begin{equation}\label{eq:IniSystem}
\Qt_a^b[\Dt_i] = I[\Dt_i], \quad i = 1,\dots, 2m
\end{equation}
and the source root $\br$ that solves (\ref{eq:IniSystem}) is the classical polynomial Gaussian rule
computed on every element.

At every instant, a certain domain $\Om \in \R^{2m}$ bounds the root.
For the source domain $\Omt \subset \R^{2m}$ we know that every element contains $\frac{d+1}{2}$ nodes and therefore
\begin{equation}\label{eq:OmegaNodes}
\begin{array}{cccc}
(\tau_1, \dots, \tau_m) \in & \underbrace{[\xt_0,\xt_1]\times \dots \times [\xt_0, \xt_{1}]}
& \times \dots \times  & \underbrace{[\xt_{n-1}, \xt_{n}] \times \dots \times [\xt_{n-1}, \xt_n].}\\
& \frac{d+1}{2} & & \frac{d+1}{2}
\end{array}
\end{equation}
For the weights we use (a rough) range $[0,b-a]$. Combined together, the source domain is
\begin{equation}\label{eq:OmegaT}
\begin{array}{cccc}
\Omt = & \underbrace{[\xt_0,\xt_1] \times \dots \times [\xt_{n-1}, \xt_n]} & \times & \underbrace{[0,b-a] \times \dots \times [0,b-a]}.\\
& m & & m
\end{array}
\end{equation}
As the source space continuously evolves to the target one, the system $\FF(\bx,t)=\mathbf 0$ continuously changes too,
and so does $\Om(t)$. The root $\br(0)$ of $\FF(\bx,0)= \mathbf 0$ is numerically traced and the root $\br(1)$
of $\FF(\bx,1)= \mathbf 0$ is returned. We refer the reader to \cite[Section 4]{Homotopy-2015} for a detailed description
of this numerical tracing.

\subsection{Limiting algebraic system}\label{sec:limit}

\begin{figure}[!tb]
\vrule width0pt\hfill
 \begin{overpic}[width=.82\columnwidth,angle=0]{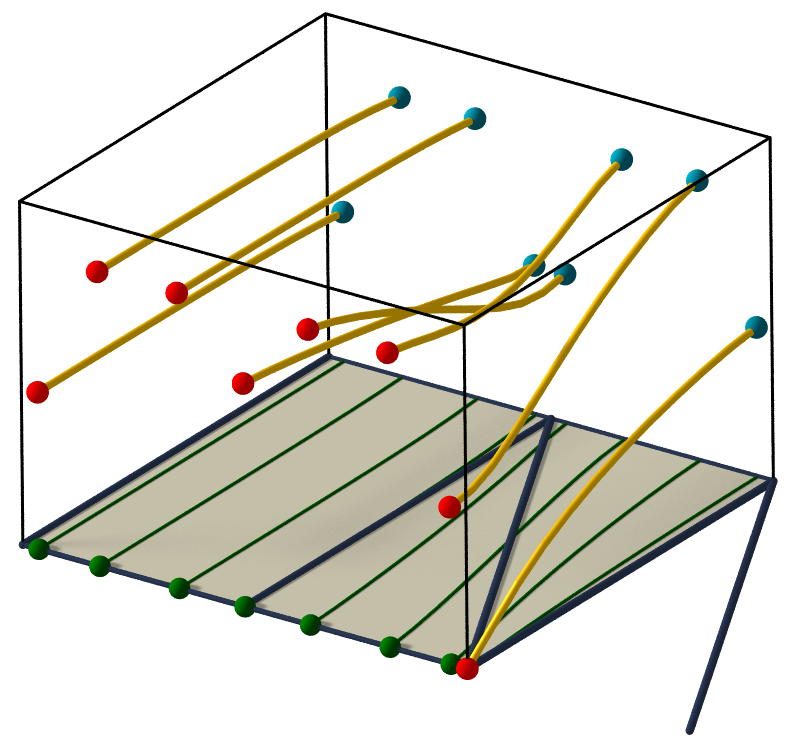}
   \put(5,5){$\bm=(8,6,8,2)$}
   \put(5,10){$\bx=(0,\frac{1}{2},1,\frac{3}{2})$}
   \put(70,91){$\xxt=(0,\frac{1}{2},1)$}
   \put(70,86){$\bmt=(8,8,8)$}
   \put(2,21){$\tau_1$}
   \put(8,19){$\tau_2$}
   \put(55,5){$b=\tau_8$}
   \put(88,1){$b+\frac{1}{N}$}
   \put(-1,24){$a$}
   \put(95,78){$[b,0,0.2]$}
   \put(-6,71){\fcolorbox{gray}{white}{\small $t=1$}}
   \put(30,93){\fcolorbox{gray}{white}{\small $t=0$}}
	\end{overpic}
\hfill \vrule width0pt\\[2ex]
\vrule width0pt\hfill
 \begin{overpic}[width=.49\columnwidth,angle=0]{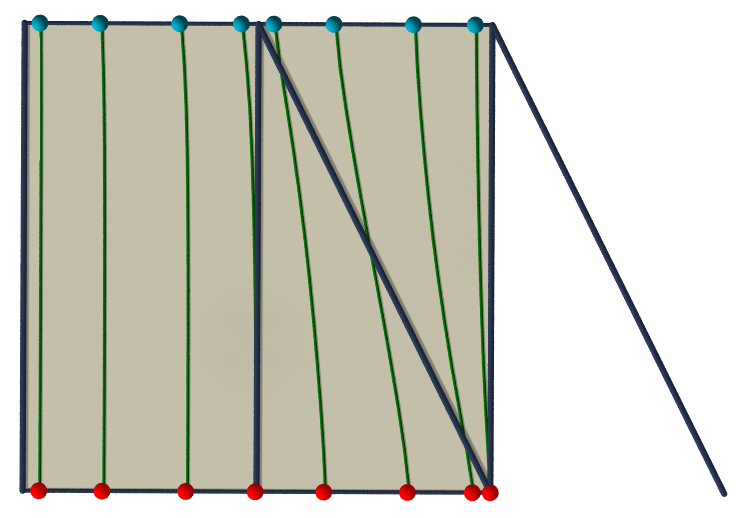}
	\end{overpic}
 \hfill
 \begin{overpic}[width=.49\columnwidth,angle=0]{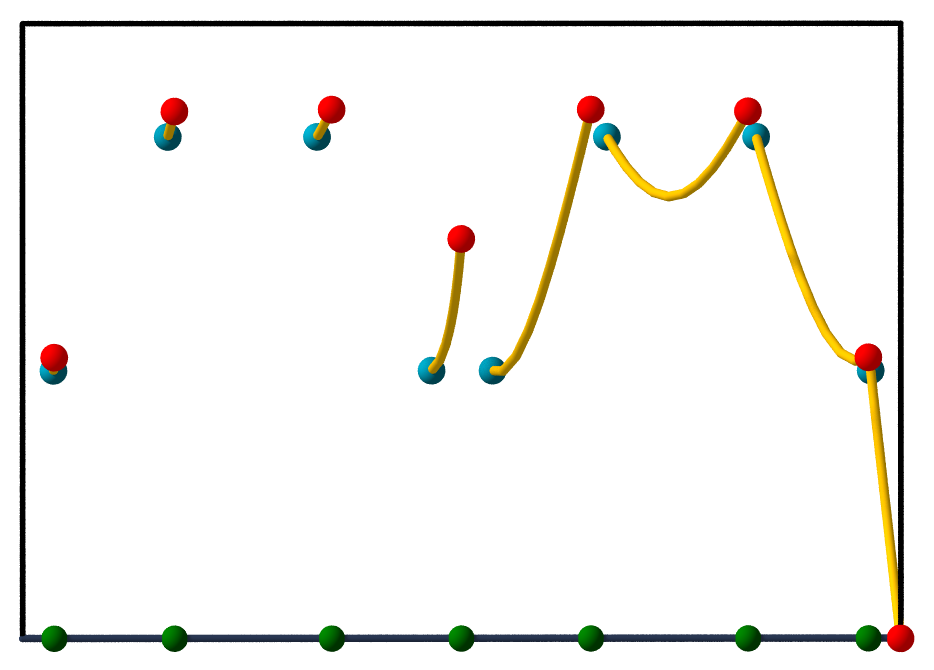}
    \put(40,50){$[\tau_4,\omega_4]$}
\end{overpic}\hfill \vrule width0pt\\[-2ex]
\Acaption{1ex}{The evolution of the Gaussian quadrature rule (\ref{quadrature}) for $d=7$, $c=1$, $N=2$, $n=2$ on the interval $[a,b]$.
The rule corresponds to the transformation of the spline spaces shown in Fig.~\ref{fig:Trans7}.
Top: starting from the optimal source rule (blue) for the space $\St^{n,d}_{\xxt,\bmt}$ with partition $\xxt$
and multiplicities $\bmt$, the evolution of the nodes (green) and weights (yellow) is shown for a geodesic knot transformation (dark blue).
Bottom left: top view of the evolution. The dimension difference, $r=2$, requires one double knot to get transformed to $b$
  while another double knot moves to $b+\frac{1}{N}=\frac{3}{2}$. Bottom right: front view of the evolution. The source polynomial Gauss rule on $n=2$ elements (blue)
  gets transformed to the target rule (red).
  At the limit stage of the continuation, $t=1$, two basis functions loose their support on $[a,b]$ which results in $\tau_{8} = b$ and $\omega_{8} = 0$.
  }\label{fig:N2}
 \end{figure}

We now discuss the behavior of the polynomial system (\ref{eq:SystemTime}) as $t\rightarrow 1$.
We already know from our homotopic setting (\ref{eq:r}) that $r$ knots and consequently
$r$ basis functions leave $[a,b]$. These $r$ functions loose their support on $[a,b]$
by moving $x_{N-1}\rightarrow b$, see Fig.~\ref{fig:Trans7}. However,
the optimal rule still has to exactly integrate these $r$ basis functions for any instant $t<1$ ($x_{N-1}<b$),
and hence $\frac{r}{2}$ nodes of the optimal rule must lie in $(x_{N-1},b)$.
Therefore, as $x_{N-1}\rightarrow b$, $\tau_i\rightarrow b$, $i=m-\frac{r}{2}+1,\dots, m$.
Additionally, we observe the following: consider, e.g. $r=2$, the last basis function $D_{2m}$ with
a non-zero support over knots
\begin{equation}\label{eq:knotsD2m}
\begin{array}{ccc}
x_{N-1}, & \underbrace{x_{N}, \dots, x_N}, & x_{N+1}, x_{N+1}\\
 & d-1 &
\end{array}
\end{equation}
and denote $h:=b-x_{N-1}$. $D_{2m}$ has a root of multiplicity $d$ at $x_{N-1}$ and therefore, by Taylor's theorem, there exists $c_1\neq 0$
such that
\begin{equation}\label{eq:D2m}
D_{2m} = c_1h^d + \mathcal O(h^{d+1})
\end{equation}
and consequently
\begin{equation}\label{eq:ID2m}
I[D_{2m}] = \frac{c_1}{d+1}h^{d+1} + \mathcal O(h^{d+2})
\end{equation}
over $[a,b]$. 
Since the rule must integrate exactly $D_{2m}$ as $h\rightarrow 0$,
we obtain from (\ref{eq:D2m}) and (\ref{eq:ID2m}) that $\omega_{m}\rightarrow 0$.

A similar logic applies to the remaining $r-1$ basis functions.
The last $r$ equations of the system (\ref{eq:SystemTime}) become ill-posed ($0=0$) as $h\rightarrow0$.
However, the above described observations, supported also by numerical results shown in Section~\ref{sec:ex}, allow us to reduce the system by these $r$ equations
to a well-constrained $(2m-r)\times (2m-r)$ system by setting
\begin{equation}\label{eq:LimitNodes}
\renewcommand{\arraystretch}{1.4}
\begin{array}{rcc}
\tau_{m-\frac{r}{2}+i}&\rightarrow b & \multirow{2}{*}{$\quad \textnormal{as} \quad t \rightarrow 1 \quad \textnormal{for} \quad i=1,\dots,\frac{r}{2}.$}\\
\omega_{m-\frac{r}{2}+i}&\rightarrow 0 &
\end{array}
\end{equation}

\section{Numerical validation of proposed Gaussian rules}\label{sec:ex}

 In this section, we show the results of our algorithm and
derive Gaussian quadrature rules for various target spaces $S^{N,d}_{\bx,\bm}$.
The evolution of the Gaussian rule for $d=7$, $c=1$, $N=2$ is shown in Fig.~\ref{fig:N2}.
The source rule starts with $m=8$ nodes while the optimal target rule requires only $7$ nodes.
As $t\rightarrow 1$, $\tau_8\rightarrow b$, i.e., this node looses its relevance on $[a,b]$.
We display the evolution of the weights as a set of space curves, showing that also $w_8\rightarrow 0$
as $t\rightarrow 1$.

The error of the rule $\Q$ is measured in terms of the Euclidean norm of the vector of the residues of the
system (\ref{eq:IniSystem}), normalized by the dimension of the system
\begin{equation}\label{eq:Error}
\|\br\| = \frac{1}{2m}(\sum_{i=1}^{2m} (\Q_a^b[D_i] - I[D_i])^2)^{\frac{1}{2}}.
\end{equation}

The optimal count of the quadrature points (\ref{eq:Micchelli}), for uniform $C^c$-continuous spaces
consisting of $N$ uniform elements, becomes
\begin{equation}\label{eq:Micchelli2}
     d + 1 + (N-1)(d-c) = 2 m.
\end{equation}
As $N\rightarrow\infty$, this count determines a certain nodal pattern.
That is, a number of optimal quadrature points required per element.
For example, for $d=5$, $c=1$, Eq.~(\ref{eq:Micchelli2}) gives a two-nodes pattern per element.
Moreover, as $N\rightarrow\infty$, we have one-element repetition for odd continuities and a two-element
repetition for even continuities.
Our rules over finite number of elements $N$ converge to their asymptotic counterparts,
thus we classify the derived rules according to these repetition patterns.

\subsection{Optimal quadrature rules for uniform, $C^1$-continuous spaces}\label{ssec:ExC1}

%
Here, the original spline space is $C^2$ continuous, i.e., $k=2$ in Section~\ref{sec:matrices}.
The optimal asymptotic rules, i.e., when $N\rightarrow\infty$, have a simple pattern: for degree $d=2s+1$, Eq.~(\ref{eq:Micchelli2}) gives
an $s$-nodes pattern per element. We now show how the derived rules for finite intervals $[0,N]$ converge to this asymptotic pattern.

$C^1$ Quintics, $d=5$, $c=1$. For uniform quintic spaces, we have recently derived a close form formula that computes recursively the optimal nodes and weights \cite{Quadrature51-2014}.
Our exact solution validates the scheme we describe herein.
Table~\ref{tab:51} displays the results derived by our homotopy continuation method
which match the analytic solution with the precision of $20$ decimal digits, cf. Fig.~16 and Eq.~39 in \cite{Quadrature51-2014}.
The nodes of the asymptotic rule are the knots and the midpoints of the elements, respectively, and their weights are
\begin{equation}\label{eq:InfyWeights51}
\omega_{2i} = \frac{7}{15} = 0.4\overline{6}, \quad \omega_{2i+1} = \frac{8}{15} = 0.5\overline{3}.
\end{equation}
In Table~\ref{tab:51}, we observe a fast convergence of the derived rule to the asymptotic counterpart since only the nodes and weights on the first four elements
differ from the asymptotic values.

\begin{table}[!tb]
 \begin{center}
  \begin{minipage}{0.9\textwidth}
\caption{Two-node patterned Gaussian quadrature rule
for $d=5$, $c=1$, with $N=10$ uniform elements over $[0,N]$.
The nodes and weights are displayed with $20$ decimal digits, showing only the nodes and weights in the first four elements differ from
the asymptotic values derived in \cite{Quadrature51-2014}, Eq.~39, by more than double machine precision.
}\label{tab:51}
  \end{minipage}
\vspace{0.2cm}\\
\small{
\renewcommand{\arraystretch}{1.15}
\renewcommand{\tf}{\small}
\begin{tabular}{| c | c || l| l|}\hline
\multicolumn{2}{|c}{} & \multicolumn{2}{c|}{$d=5$, $c=1$, $N=10$, uniform, $\|\br\| = 2.37^{-18}$ } \\\hline
$\#$el. & $i$ &  \multicolumn{1}{c|}{$\tau_i$}  & \multicolumn{1}{c|}{$\omega_i$} \\\hline\hline
\multirow{2}{*}{1} & 1 & \tf 0.12251482265544137787 & \tf 0.30201742881457235729   \\
                   & 2 & \tf 0.54415184401122528880 & \tf 0.48501960822246467975   \\\hline
\multirow{2}{*}{2} & 3 & \tf 1.00646547160565963977 & \tf 0.44671772013629118653   \\
                   & 4 & \tf 1.50027307286873389123 & \tf 0.53303872093804185483   \\\hline
\multirow{2}{*}{3} & 5 & \tf 2.00003879729563044051 & \tf 0.46653987137191212073   \\
                   & 6 & \tf 2.50000001053211375767 & \tf 0.53333332209820754959   \\\hline
\multirow{2}{*}{4} & 7 & \tf 3.00000000150452933969 & \tf 0.46666666175184358463   \\
                   & 8 & \tf 3.5                    & \tf 0.53333333333333333333   \\\hline
\multirow{2}{*}{5} & 9 & \tf 4                      & \tf 0.46666666666666666667   \\
                   & 10& \tf 4.5                    & \tf 0.53333333333333333333   \\\hline
\end{tabular}
}
\end{center}
\end{table}

\begin{figure}[!tb]
 \vrule width0pt\hfill
 \begin{overpic}[width=.99\columnwidth,angle=0]{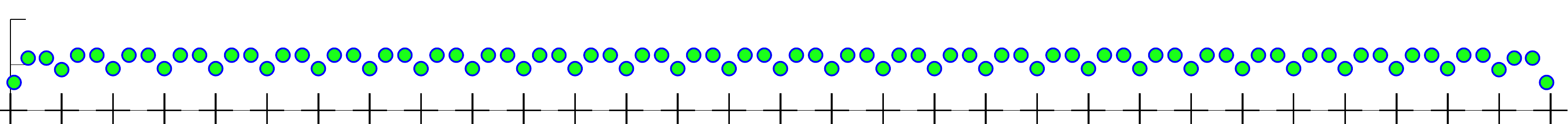}
    \put(20,9){\fcolorbox{gray}{white}{$d=7$, $c=1$, $N=30$, $\card(\XX_N) = 190$, $r=2$}}
    \put(3,6){\footnotesize $0.6$}
    \put(0,-2){\footnotesize $0$}
    \put(90,-2){ $x_{30}=30$}
	\end{overpic}
 \hfill \vrule width0pt\\[3ex]
 \vrule width0pt\hfill
 \begin{overpic}[width=.89\columnwidth,angle=0]{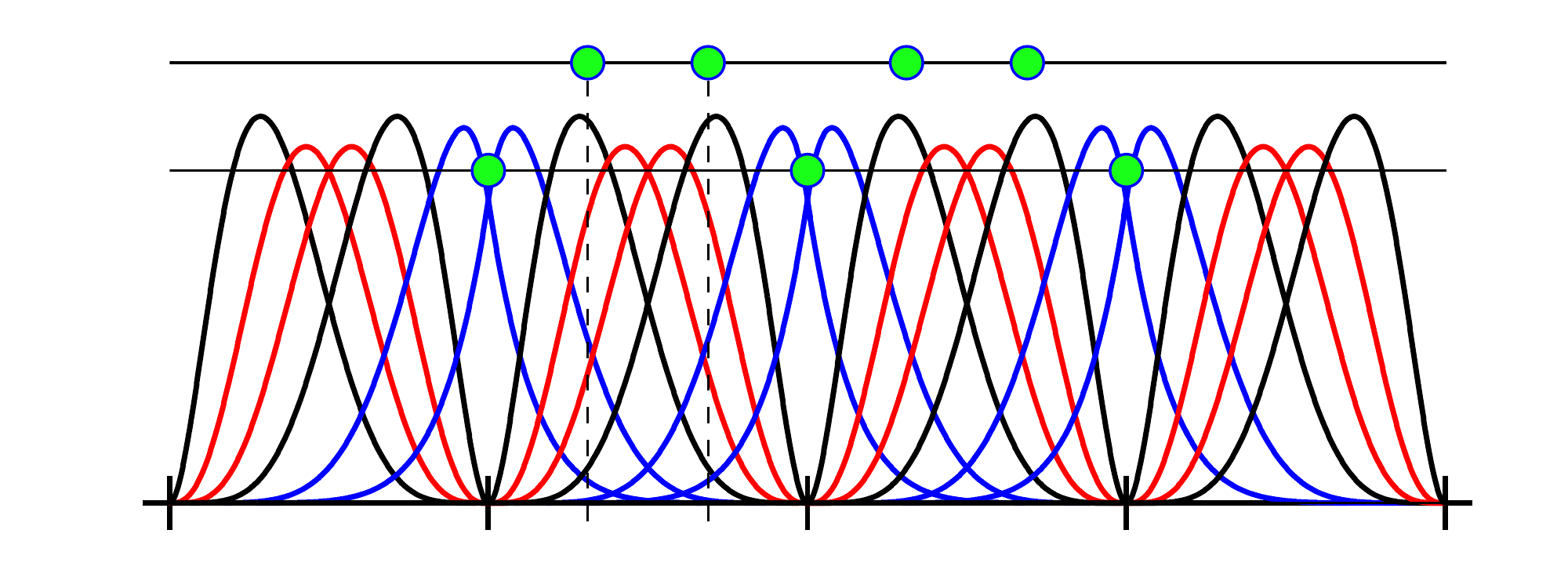}
    \put(8,1){$x_{i-1}$}
    \put(29,1){$x_i$}
    \put(33,-1){$d_1$}
    \put(47,-1){$d_1$}
    \put(51,1){$x_{i+1}$}
    \put(70,1){$x_{i+2}$}
    \put(31,4){$\scriptsize\underbrace{\makebox[0.7cm]{}}$}
    \put(45,4){$\scriptsize\underbrace{\makebox[0.7cm]{}}$}
    \put(95,33){$\frac{49}{135}$}
    \put(95,26){$\frac{37}{135}$}
	\end{overpic}
 \hfill \vrule width0pt\\[-4ex]
 \Acaption{1ex}{Three-node-per-element rule ($d=7$, $c=1$).
 Top: the layout of the optimal quadrature rule (green dots)
 for uniform knot distribution with $N=30$ elements over $[a,b]=[0,30]$ is shown, see also Table~\ref{tab:71}.
 Bottom: Asymptotic layout of the optimal rule.
 Two types of basis functions have support only on one element (red and black), whilst another type spans two (blue).
 The arrangement of the basis functions is symmetric with respect to the elements' midpoints and knots.
 One set of nodes are the knots
 with the weight $\frac{37}{135}$. The position of the nodes lying inside the elements is determined by $d_1$.
 The system for the interior nodes with $x_i=0$ and $x_{i+1}=1$ is shown in (\ref{eq:LimSys71}).
 \vspace{0.5cm}
 }\label{fig:71}
 \end{figure}

 $C^1$ Septics, $d=7$, $c=1$. For uniform septic spline space, our optimal quadrature rule with $N=30$ elements is shown in Fig.~\ref{fig:71} top,
 see also Table~\ref{tab:71}. For this space there is not explicit Gaussian rule to validate our method against.
 We see again that only nodes and weights on a few boundary elements differ from a regular pattern formed inside the interval. 
 In the asymptotic setting, i.e., when integrating over a whole real line ($N=\infty$),
 Eq.~(\ref{eq:Micchelli2}) gives us that three nodes per element are required.
 We build the asymptotic system similarly to \cite{Hughes-2010}.
 Due to the symmetry of the basis functions,
 one set of nodes are the knots and the other lie inside the elements, arranged symmetrically with respect to the elements' midpoints,
 see Fig.~\ref{fig:71} bottom. The position of the nodes lying inside the elements is determined by $d_1$
 and, for the particular (normalized) setting $x_i=0$ and $x_{i+1}=1$, we obtain the asymptotic system for the interior node and weight as
\begin{equation}\label{eq:LimSys71}
    \renewcommand{\arraystretch}{1.4}
    \begin{array}{ccc}
        21w_1d_1^2(1-d_1)^5+21w_1(1-d_1)^2d_1^5 & = & \frac{1}{8},  \\
        35w_1d_1^3(1-d_1)^4+35w_1(1-d_1)^3 d_1^4& = & \frac{1}{8},
    \end{array}
\end{equation}
 which corresponds to the $\Q[D_i] = I[D_i]$ constraint for two basis functions (red and black in Fig.~\ref{fig:71} bottom).
Solving (\ref{eq:LimSys71}), we obtain
 \begin{equation}\label{eq:InfyNodes}
 d_1 = \frac{7-\sqrt{7}}{14} \doteq 0.31101776349538638639
 \end{equation}
 and
  \begin{equation}\label{eq:InfyWeights}
 \renewcommand{\arraystretch}{1.45}
  \begin{array}{rcr}
\omega_{4i} = \omega_{4i+3} & = & \frac{37}{135} = 0.27\overline{407}, \\
\omega_{4i+1} = \omega_{4i+2} & = & \frac{49}{135} = 0.36\overline{296}. \\
\end{array}
 \end{equation}
 Therefore the asymptotic Gaussian rule ($N=\infty$) is
\begin{equation}\label{eq:limit71}
\int_{\mathbb R}f(t)\, \mathrm{d}t = \sum_{i \in \mathbb Z} h(\frac{37}{135} f(ih) + \frac{49}{135} f(h(i+d_1)) + \frac{49}{135} f(h(i+1-d_1))),
\end{equation}
where $h$ is the size of the element and $d_1$ is the limit value from (\ref{eq:InfyNodes}). 
 From Table~\ref{tab:71}, we observe a fast convergence of our rule computed for $N=30$
 to its asymptotic counterpart (\ref{eq:limit71}). With the precision of $20$ decimal digits,
 only the nodes and weights on the first four elements differ from the asymptotic values.
 
 \begin{table}[!tb]
\begin{center}
\begin{minipage}{0.9\textwidth}
    \caption{Nodes and weights for Gaussian quadrature rules (\ref{quadrature}) for $C^1$ septics, $d=7$, $c=1$, and $N=30$ uniform elements over $[0,30]$,
    shown in Fig.~\ref{fig:71} top.
    The nodes and weights from the first six boundary elements are displayed with $20$ decimal digits, showing the convergence
    to the asymptotic rule (\ref{eq:limit71}), cf. (\ref{eq:InfyNodes}) and (\ref{eq:InfyWeights}).}\label{tab:71}
\end{minipage}
\vspace{0.2cm}\\
\small{
\renewcommand{\arraystretch}{1.15}
\renewcommand{\tf}{\small}
\begin{tabular}{| c | c || l| l|}\hline
\multicolumn{2}{|c}{} & \multicolumn{2}{c|}{$d=7$, $c=1$, $N=30$, uniform, $\|\br\| = 1.09^{-18}$ } \\\hline
$\#$el. & $i$ &  \multicolumn{1}{c|}{$\tau_i$}  & \multicolumn{1}{c|}{$\omega_i$} \\\hline\hline
\multirow{3}{*}{1} & 1 & \tf 0.07299402407314973216 & \tf 0.18285701415655202878  \\
                   & 2 & \tf 0.34700376603835188472 & \tf 0.34297577246926732566  \\
                   & 3 & \tf 0.70500220988849838312 & \tf 0.34416721337418064556  \\\hline
\multirow{3}{*}{2} & 4 & \tf 1.00213067803177481153 & \tf 0.26713002701651926831  \\
                   & 5 & \tf 1.31109168439816575861 & \tf 0.36292347046348192394  \\
                   & 6 & \tf 1.68901548923246352193 & \tf 0.36292410619137875755  \\\hline
\multirow{3}{*}{3} & 7 & \tf 2.00000433077293358133 & \tf 0.27405943376486496347  \\
                   & 8 & \tf 2.31101776381410751148 & \tf 0.36296296279505220171  \\
                   & 9 & \tf 2.68898223664848840334 & \tf 0.36296296279505818626  \\\hline
\multirow{3}{*}{4} & 10& \tf 3.00000000001875375310 & \tf 0.27407407401068173581  \\
                   & 11& \tf 3.31101776349538638640 & \tf 0.36296296296296296296  \\
                   & 12& \tf 3.68898223650461361361 & \tf 0.36296296296296296296  \\\hline
\multirow{3}{*}{5} & 13& \tf 4                      & \tf 0.27407407407407407407  \\
                   & 14& \tf 4.31101776349538638639 & \tf 0.36296296296296296296  \\
                   & 15& \tf 4.68898223650461361361 & \tf 0.36296296296296296296  \\\hline
\multirow{3}{*}{6} & 16& \tf 5                      & \tf 0.27407407407407407407  \\
                   & 17& \tf 5.31101776349538638639 & \tf 0.36296296296296296296  \\
                   & 18& \tf 5.68898223650461361361 & \tf 0.36296296296296296296  \\\hline
\end{tabular}
}
\end{center}
\end{table}

 This limit rule requires only three Gaussian points per element whilst the classical polynomial
 Gauss (used as our source rule) requires four nodes per element. 
%
From Table~\ref{tab:71}, we also conclude that for spaces  $S^{N,7}_{\bx,\bm}$ with uniform knot vectors and large number of elements ($N>30$)
one can use our values from Table~\ref{tab:71} for the first five elements and for all the remaining interior elements use
the asymptotic rule (\ref{eq:limit71}).

\begin{figure}[!tb]
\vrule width0pt\hfill
 \begin{overpic}[width=.49\columnwidth,angle=0]{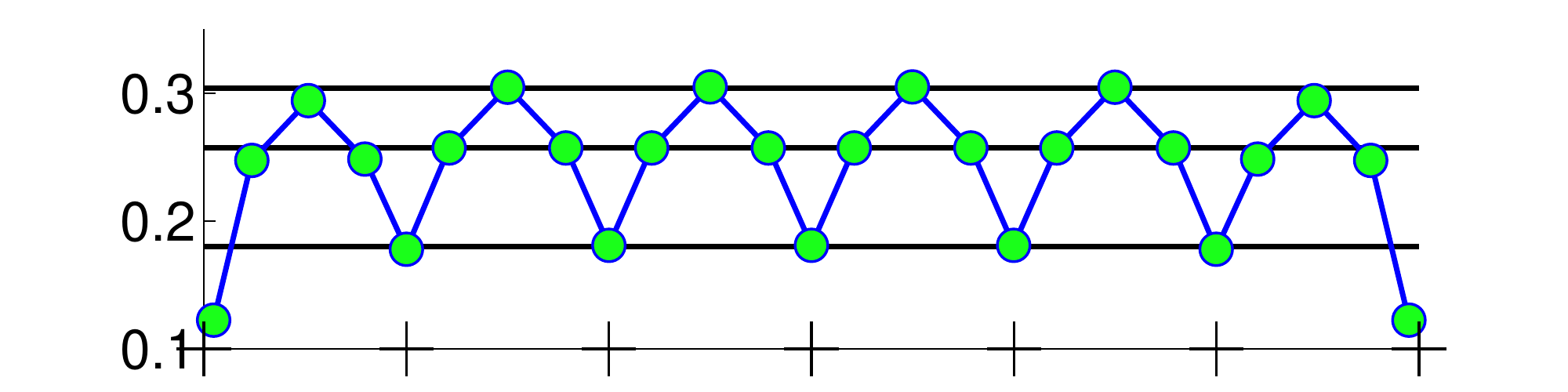}
    \put(90,25){\fcolorbox{gray}{white}{$d=9$, $c=1$}}
    \put(15,25){\fcolorbox{gray}{white}{$N=6$, $r=0$}}
    \put(82,-3){$x_{6}=6$}
	\end{overpic}
 \hfill
 \begin{overpic}[width=.49\columnwidth,angle=0]{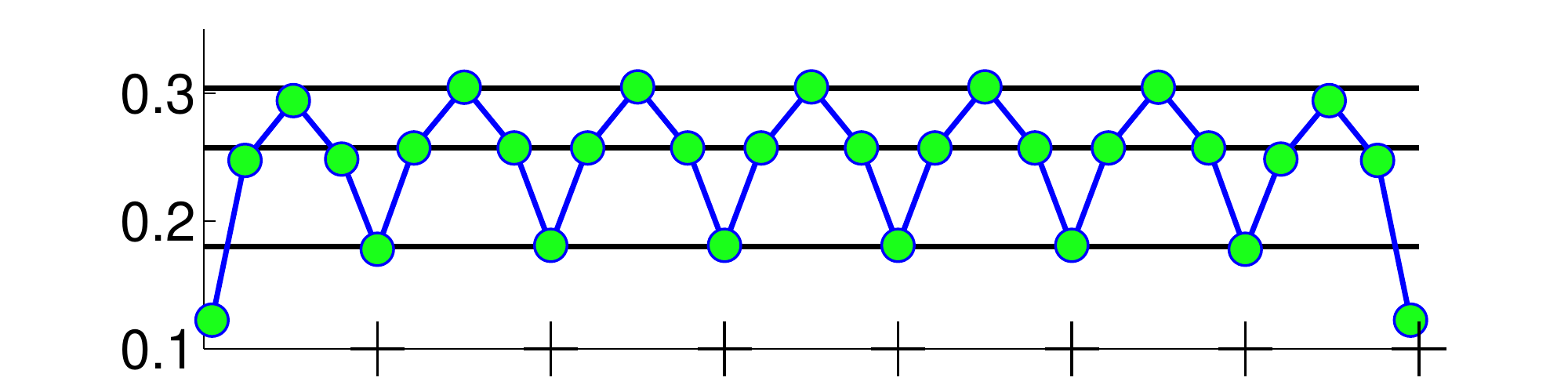}
    \put(60,25){\fcolorbox{gray}{white}{$N=7$, $r=2$}}
    \put(80,-3){$x_{7}=7$}
	\end{overpic}
 \hfill \vrule width0pt\\[-1ex]
 \vrule width0pt\hfill
 \begin{overpic}[width=.89\columnwidth,angle=0]{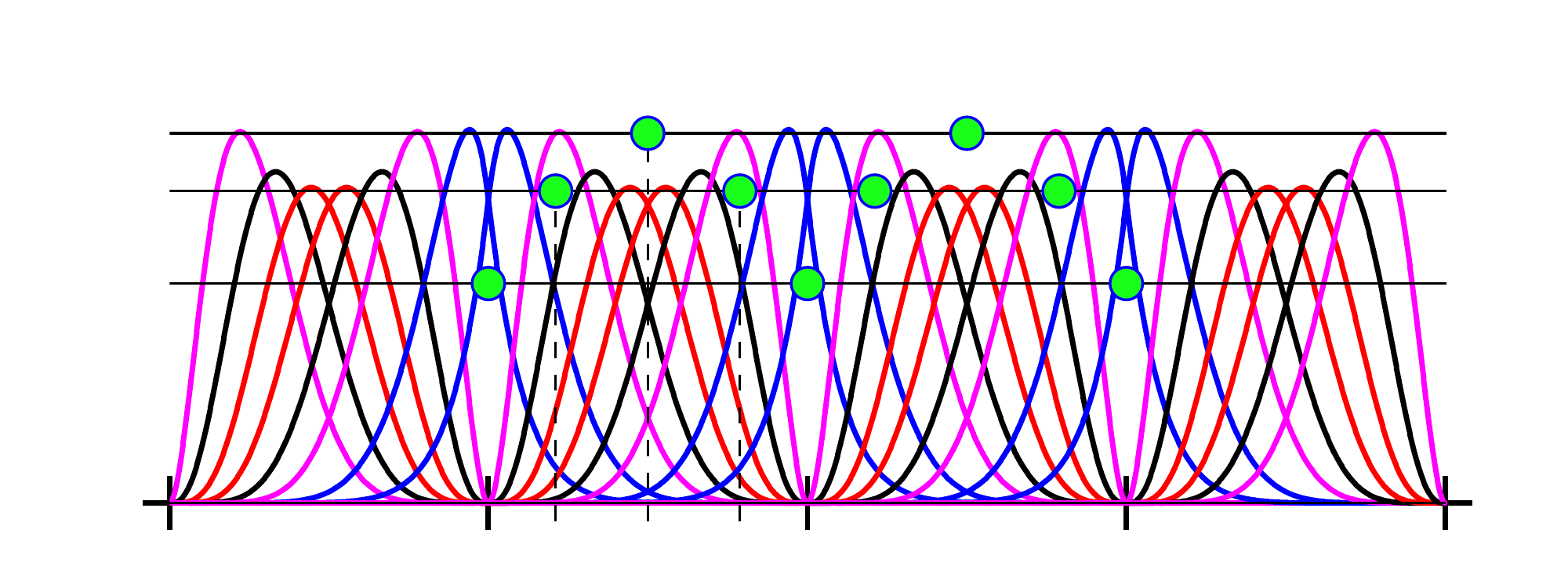}
    \put(8,1){$x_{i-1}$}
    \put(29,1){$x_i$}
    \put(33,-1){$d_1$}
    \put(47,-1){$d_1$}
    \put(51,1){$x_{i+1}$}
    \put(70,1){$x_{i+2}$}
    \put(31,4){$\scriptsize\underbrace{\makebox[0.5cm]{}}$}
    \put(47,4){$\scriptsize\underbrace{\makebox[0.5cm]{}}$}
    \put(94,29){\footnotesize $\frac{32}{105}$}
    \put(94,24){\footnotesize $\frac{27}{105}$}
    \put(94,19){\footnotesize $\frac{19}{105}$}
	\end{overpic}
 \hfill \vrule width0pt\\[-2ex]
 \Acaption{1ex}{Four-node patterned rules. Top: optimal rules for $C^1$ nonics, $d=9$, $c=1$ over finite domains $[0,N]$
 for $N=6,7$ are shown.
 The layout of the optimal quadrature rules (green dots) converges to the asymptotic pattern (bottom).
 In this limit case ($N=\infty$), one set of nodes becomes the knots, another becomes the elements' midpoints, and the last set
 of nodes lies inside, determined by the distance $d_1$, see (\ref{eq:limit91}).}\label{fig:d=9_k=1}
 \end{figure}

$C^1$ Nonics, $d=9$, $c=1$. Examples of optimal rules for nonic splines with uniform knot vectors for various number of elements are shown in Fig.~\ref{fig:d=9_k=1}.
We again observe the convergence to the asymptotic four-point rule
\begin{equation}\label{eq:limit91}
\int_{\mathbb R}f(t)\, \mathrm{d}t = \sum_{i \in \mathbb Z} h(\frac{19}{105} f(ih) + \frac{32}{105} f(i+\frac{h}{2}) + \frac{27}{105} (f(h(i+d_1)) + f(h(i+1-d_1)))).
\end{equation}
where
\begin{equation}\label{eq:limit91d}
d_1 = 0.21132486540518711775
\end{equation}
was computed numerically with $N=20$ elements.

 \begin{table}[!tb]
\begin{center}
\begin{minipage}{0.9\textwidth}
    \caption{Nodes and weights for Gaussian quadrature rules (\ref{quadrature}) for $d=9$, $c=1$ and $N=20$ uniform elements over $[0,20]$.
    The nodes and weights are displayed with $20$ decimal digits, showing the convergence
    to the asymptotic rule (\ref{eq:limit91}).
    The asymptotic weights are $\frac{19}{105} = 0.1\overline{809523}$,
    $\frac{27}{105} = 0.2\overline{571428}$, and $\frac{32}{105} = 0.3\overline{047619}$.
    }\label{tab:91}
\end{minipage}
\vspace{0.2cm}\\
\small{
\renewcommand{\arraystretch}{1.15}
\renewcommand{\tf}{\small}
\begin{tabular}{| c | c || l| l|}\hline
\multicolumn{2}{|c}{} & \multicolumn{2}{c|}{$d=9$, $c=1$, $N=20$, uniform, $\|\br\| = 3.74^{-16}$ } \\\hline
$\#$el. & $i$ &  \multicolumn{1}{c|}{$\tau_i$}  & \multicolumn{1}{c|}{$\omega_i$} \\\hline\hline
\multirow{4}{*}{1} & 1 & \tf 0.04850054944699732930  & \tf 0.12248110464981389735  \\
                   & 2 & \tf 0.23860073755186230506  & \tf 0.24745843345844748980  \\
                   & 3 & \tf 0.51704729510436750234  & \tf 0.29425875345698032366  \\
                   & 4 & \tf 0.79585141789677286330  & \tf 0.24839430102735088178  \\\hline
\multirow{4}{*}{2} & 5 & \tf 1.00090689188286066327  & \tf 0.17791129711197193832  \\
                   & 6 & \tf 1.21134881170777431791  & \tf 0.25713493030888296163  \\
                   & 7 & \tf 1.50001514554780606887  & \tf 0.30475265725234796893  \\
                   & 8 & \tf 1.78868153181393990233  & \tf 0.25713507064149301837  \\\hline
\multirow{4}{*}{3} & 9 & \tf 2.00000079645799274743  & \tf 0.18094964259082567701  \\
                   & 10& \tf 2.21132486542419774005  & \tf 0.25714285713665875486  \\
                   & 11& \tf 2.50000000001205222537  & \tf 0.30476190475455863188  \\
                   & 12& \tf 2.78867513459990674909  & \tf 0.25714285713665885746  \\\hline
\multirow{4}{*}{4} & 13& \tf 3.00000000000063432715  & \tf 0.18095238095020007515  \\
                   & 14& \tf 3.21132486540518711775  & \tf 0.25714285714285714286  \\
                   & 15& \tf 3.5                     & \tf 0.30476190476190476190  \\
                   & 16& \tf 3.78867513459481288225  & \tf 0.25714285714285714286  \\\hline
\multirow{4}{*}{5} & 17& \tf 3.00000000000000000236  & \tf 0.18095238095238095238  \\
                   & 18& \tf 3.21132486540518711775  & \tf 0.25714285714285714286  \\
                   & 19& \tf 3.5                     & \tf 0.30476190476190476190  \\
                   & 20& \tf 3.78867513459481288225  & \tf 0.25714285714285714286  \\\hline
\end{tabular}
}
\end{center}
\end{table}

\subsection{Optimal quadrature rules for uniform $C^0$-continuous spaces}\label{ssec:ExEven}

For these spaces, the asymptotic configuration forms a two-element repetition.
That is, the count of optimal nodes (\ref{eq:Micchelli2}) for $N=\infty$ becomes ``$\textnormal{odd} = 2m$''.
Therefore the number of nodes per element is not constant as in Section~\ref{ssec:ExC1}, but repetitively changes
on odd and even elements.

\begin{figure}[!tb]
\vrule width0pt\hfill
 \begin{overpic}[width=.49\columnwidth,angle=0]{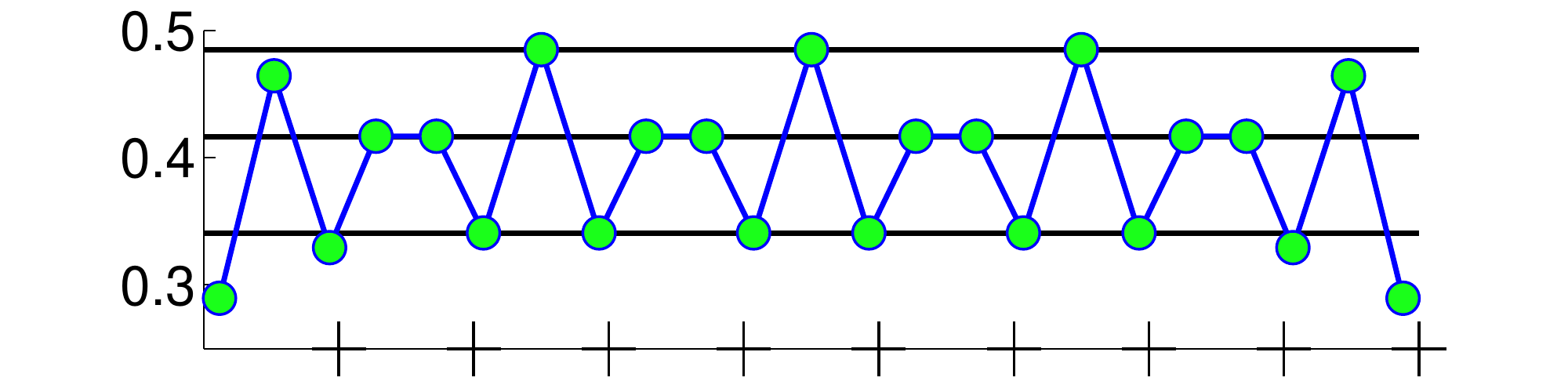}
    \put(90,30){\fcolorbox{gray}{white}{$d=5$, $c=0$}}
    \put(15,29){\fcolorbox{gray}{white}{$N=9$, $r=1$}}
    \put(82,-3){$x_{9}=9$}
	\end{overpic}
 \hfill
 \begin{overpic}[width=.49\columnwidth,angle=0]{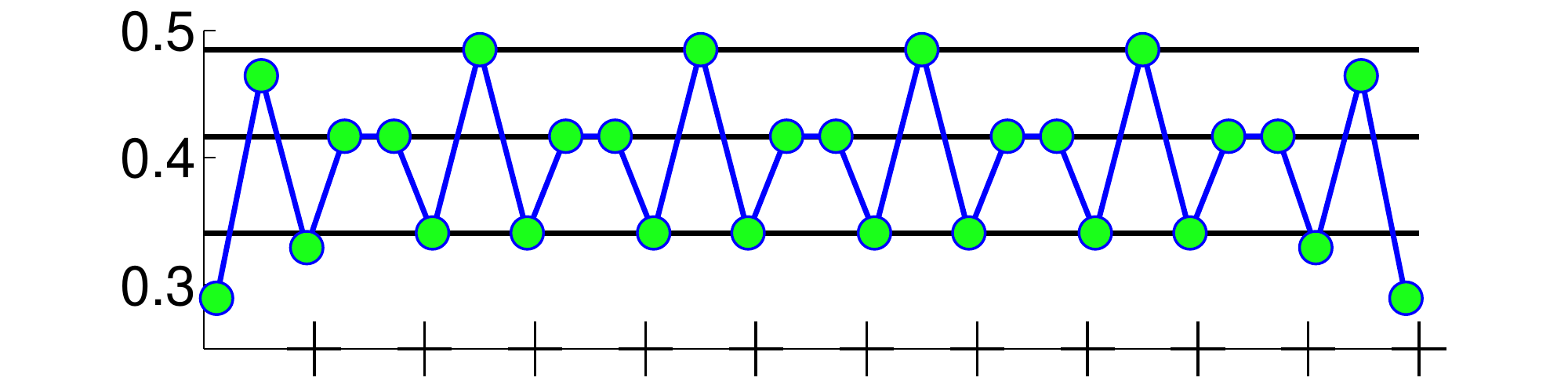}
    \put(65,29){\fcolorbox{gray}{white}{$N=11$, $r=2$}}
    \put(80,-3){$x_{11}=11$}
	\end{overpic}
 \hfill \vrule width0pt\\[2ex]
 \vrule width0pt\hfill
 \begin{overpic}[width=.89\columnwidth,angle=0]{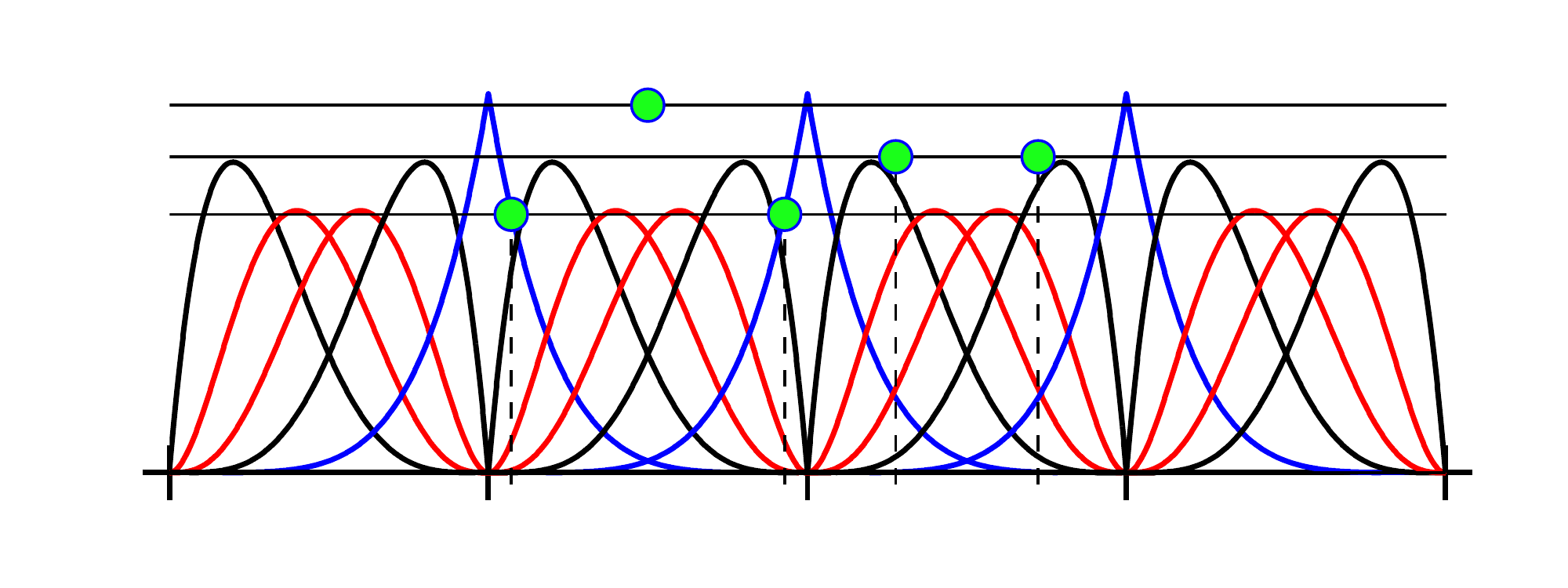}
    \put(8,3){$x_{i-1}$}
    \put(28,3){$x_i$}
    \put(53.5,1){$d_2$}
    \put(68,1){$d_2$}
    \put(46,3){$x_{i+1}$}
    \put(72,3){$x_{i+2}$}
    \put(52,6){$\scriptsize\underbrace{\makebox[0.6cm]{}}$}
    \put(66,6){$\scriptsize\underbrace{\makebox[0.6cm]{}}$}
    \put(31,5){\scalebox{0.4}{$\scriptsize\underbrace{\makebox[0.05cm]{}}$}}
    \put(31.5,1){$d_1$}
    %
    \put(5,30){\footnotesize $\frac{64}{132}$}
    \put(94,27){\footnotesize $\frac{55}{132}$}
    \put(5,23){\footnotesize $\frac{45}{132}$}
	\end{overpic}
 \hfill \vrule width0pt\\[-4ex]
 \Acaption{1ex}{Top: optimal quadrature rules (green dots) for $C^0$-continuous quintic spline spaces for various $N$ are shown.
  The rules convergence to the asymptotic layout which requires $2.5$ nodes per element. The horizontal lines are the asymptotic weights.
  Bottom: the asymptotic ($N=\infty$) configuration with five nodes per two elements. On odd elements, the three nodes are the middle of the elements
  and the two nodes inside, determined by the distance $d_1$. On even elements, the two nodes are determined by $d_2$, see (\ref{eq:limit50d}).
   }\label{fig:d=5_k=0}
 \end{figure}

$C^0$ Quintics, $d=5$, $c=0$.
Analogously to the case for $d=7$, $c=1$, see Eq.~(\ref{eq:LimSys71}), one can build an algebraic system of equations for the asymptotic configuration.
Solving such a system, the asymptotic $2.5$-nodes-per-element rule is
\begin{equation}\label{eq:limit50}
\begin{split}
\int_{\mathbb R}f(t) & = \sum\limits_{i \in \mathbb Z} h( \frac{45}{132} (f((2i+d_1)h) + f((2i+1-d_1)h)) + \frac{64}{132} f(\frac{4i+1}{2}h)  \\
 & + \frac{55}{132} (f(h(2i+1+d_2)) + f(h(2i+2-d_2))))
\end{split}
\end{equation}
where
\begin{equation}\label{eq:limit50d}
\renewcommand{\arraystretch}{1.3}
\begin{array}{lcl}
d_1 & = & 0.07182558071116236600,  \\
d_2 & = & 0.27639320225002103036.
\end{array}
\end{equation}
The asymptotic rule is shown in Fig.~\ref{fig:d=5_k=0} bottom. A fast convergence of the rules for finite domains to the asymptotic counterpart
is shown in Table~\ref{tab:50} where for $N=11$ only the nodes and weights on the first element differ from the asymptotic values.
We conclude that for $N>11$, the optimal rule at hand consists of the first element block (lines 1 to 3) in Table~\ref{tab:50}
and the asymptotic $2.5$-nodes-per-element rule (\ref{eq:limit50}).

\begin{table}[!tb]
\begin{center}
\begin{minipage}{0.9\textwidth}
\caption{$2.5$-nodes-per-element Gaussian quadrature rule.
The nodes and weights for $d=5$, $c=0$, with $N=11$ uniform elements over $[0,N]$
are displayed with $20$ decimal digits, showing only the nodes and weights in the first one element differ from
the values of the asymptotic rule (\ref{eq:limit50}). The asymptotic weights are $\frac{45}{132} = 0.34\overline{09}$,
$\frac{55}{132} = 0.41\overline{6}$, and $\frac{64}{132} = 0.\overline{48}$.
}\label{tab:50}
  \end{minipage}
\vspace{0.2cm}\\
\small{
\renewcommand{\arraystretch}{1.15}
\renewcommand{\tf}{\small}
\begin{tabular}{| c | c || l| l|}\hline
\multicolumn{2}{|c}{} & \multicolumn{2}{c|}{$d=5$, $c=0$, $N=11$, uniform, $\|\br\| = 5.44^{-24}$ } \\\hline
$\#$el. & $i$ &  \multicolumn{1}{c|}{$\tau_i$}  & \multicolumn{1}{c|}{$\omega_i$} \\\hline\hline
\multirow{3}{*}{1} & 1 & \tf 0.11750568765381614665 & \tf 0.28964324113876964987    \\
                   & 2 & \tf 0.52159764492771382263 & \tf 0.46424613902185656449    \\
                   & 3 & \tf 0.93232523884704145929 & \tf 0.32944395317270711897    \\\hline
\multirow{2}{*}{2} & 4 & \tf 1.27639320225002103036 & \tf 0.41666666666666666667    \\
                   & 5 & \tf 1.72360679774997896964 & \tf 0.41666666666666666667    \\\hline
\multirow{3}{*}{3} & 6 & \tf 2.07182558071116236600 & \tf 0.34090909090909090909    \\
                   & 7 & \tf 2.5                    & \tf 0.48484848484848484848    \\
                   & 8 & \tf 2.92817441928883763400 & \tf 0.34090909090909090909    \\\hline
\multirow{2}{*}{4} & 9 & \tf 3.27639320225002103036 & \tf 0.41666666666666666667    \\
                   & 10& \tf 3.72360679774997896964 & \tf 0.41666666666666666667    \\\hline
\multirow{3}{*}{5} & 11& \tf 4.07182558071116236600 & \tf 0.34090909090909090909    \\
                   & 12& \tf 4.5                    & \tf 0.48484848484848484848    \\
                   & 13& \tf 4.92817441928883763400 & \tf 0.34090909090909090909    \\\hline
\end{tabular}
}
\end{center}
\end{table}

\begin{figure}[!tb]
 \vrule width0pt\hfill
 \begin{overpic}[width=.49\columnwidth,angle=0]{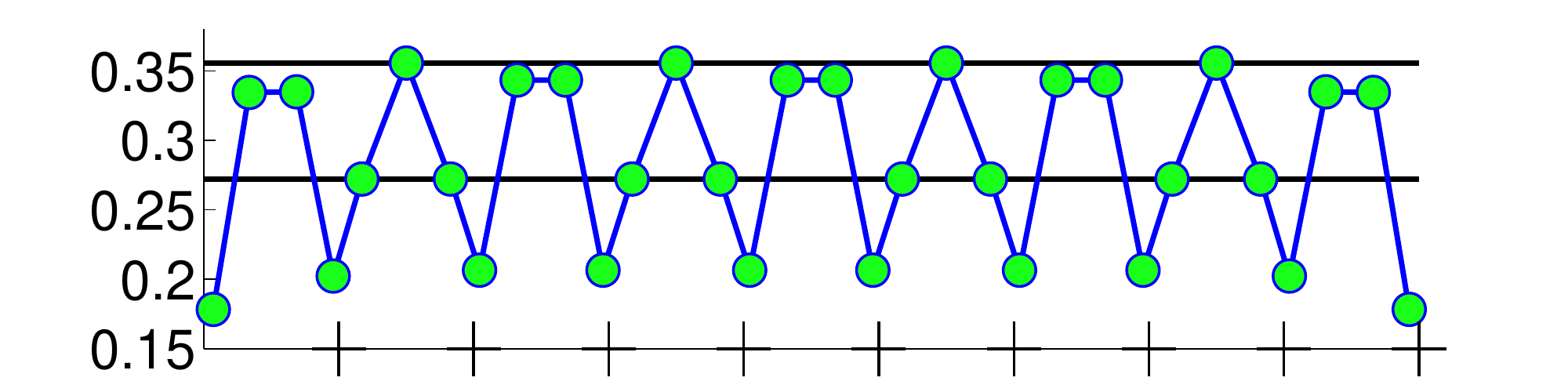}
    \put(90,30){\fcolorbox{gray}{white}{$d=7$, $c=0$}}
    \put(15,28){\fcolorbox{gray}{white}{$N=9$, $r=0$}}
    \put(82,-3){$x_{9}=9$}
	\end{overpic}
 \hfill
 \begin{overpic}[width=.49\columnwidth,angle=0]{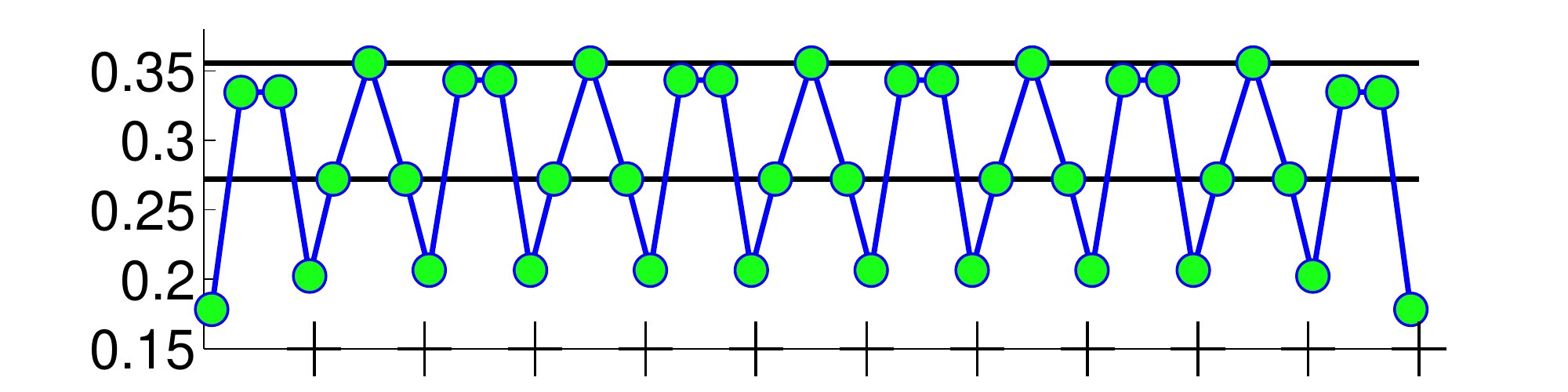}
    \put(65,28){\fcolorbox{gray}{white}{$N=11$, $r=2$}}
    \put(80,-3){$x_{11}=11$}
    \put(91,19){\small $\frac{64}{180}$}
    \put(91,10){\small $\frac{49}{180}$}
	\end{overpic}
 \hfill \vrule width0pt\\[-1ex]
 \Acaption{1ex}{Optimal quadrature rules (green dots) for septic $C^0$ splines ($d=7$, $c=0$) with uniform knots for various number of elements $N$
  on interval $[0,N]$ are shown. The asymptotic number of nodes per element is $3.5$, see (\ref{eq:Micchelli2}), and for $N=11$ the layout of nodes
  already approaches the asymptotic pattern, see also Table~\ref{tab:70},
  repetitively changing from three to four.}\label{fig:d=7_k=0}
 \end{figure}

$C^0$ Septics, $d=7$, $c=0$.
Fig.~\ref{fig:d=7_k=0} shows optimal rules for $C^0$-continuous septic splines with uniform knots for various
number of elements $N$.
Analogously to (\ref{eq:limit71}), one can derive the asymptotic rule for this case.
The asymptotic pattern is $3.5$ nodes per element.
For the three-nodes elements, the asymptotic system (over a unit element) gives
\begin{equation*}
\renewcommand{\arraystretch}{1.4}
\begin{array}{cclcccl}
d_1  & = & \frac{1}{2} - \frac{\sqrt{21}}{14},  & \quad & w_1 & = &\frac{49}{180},\\
d_2  & = & \frac{1}{2},                         & \quad & w_2 & = &\frac{64}{180},\\
d_3  & = & \frac{1}{2} + \frac{\sqrt{21}}{14},  & \quad & w_3 & = &\frac{49}{180}. \\
\end{array}
\end{equation*}
For the four-nodes elements, one builds a $4\times4$ system (2 nodes and 2 weights as the unknowns)
and, using a computer algebra software to factor the system, the problem reduces to a univariate polynomial
\begin{equation}\label{eq:quartic}
112 x^4 - 224 x^3 + 141 x^2 - 29x + 1,
\end{equation}
where its roots are the positions of the nodes on the unit interval.
Computing these roots numerically, we obtain
\begin{equation*}
\renewcommand{\arraystretch}{1.1}
\begin{array}{ccc}
d_1  & = & 0.04279465186386840500,   \\
d_2  & = & 0.32101760363894084659,   \\
d_3  & = & 0.67898239636105915341,   \\
d_4  & = & 0.95720534813613159500,
\end{array}
\end{equation*}
which match our results with 20 decimal digits, see Table~\ref{tab:70}.
The corresponding weights can be computed accordingly, however,
their explicit formula (compared to (\ref{eq:quartic})) is tedious and we skip it for the sake of brevity.
From Table~\ref{tab:70}, we again observe a very fast convergence to the asymptotic configuration ($N=\infty$)
and, with precision of 20 decimal digits, only the nodes and weights in
the first element (lines $1$ to $4$ in Table~\ref{tab:70}) differ from the asymptotic values.
We conclude that for $N>11$, the optimal rule is given by the values for the first element
and the repetitive $3.5$ pattern (lines $5$ to $11$ in Table~\ref{tab:70}).

\begin{table}[!tb]
 \begin{center}
  \begin{minipage}{0.9\textwidth}
\caption{Nodes and weights for Gaussian quadrature rule (\ref{quadrature}) for $C^0$ septics, $d=7$, $c=0$, with uniform and normalized $N=11$ elements,
shown in Fig.~\ref{fig:d=7_k=0} bottom.
The nodes and weights are displayed with $20$ decimal digits, showing only the nodes and weights in the first element differ from
the asymptotic values.}\label{tab:70}
  \end{minipage}
\vspace{0.2cm}\\
\small{
\renewcommand{\arraystretch}{1.15}
\renewcommand{\tf}{\small}
\begin{tabular}{| c | c || l| l|}\hline
\multicolumn{2}{|c}{} & \multicolumn{2}{c|}{$d=7$, $c=0$, $N=11$, uniform, $\|\br\| = 6.23^{-21}$ } \\\hline
$\#$el. & $i$ &  \multicolumn{1}{c|}{$\tau_i$}  & \multicolumn{1}{c|}{$\omega_i$} \\\hline\hline
\multirow{4}{*}{1} & 1 & \tf 0.07119064594392293204 & \tf 0.17833627514522273858   \\
                   & 2 & \tf 0.33839943987188877247 & \tf 0.33441583955888286197   \\
                   & 3 & \tf 0.68726340145291639294 & \tf 0.33496578542130575465   \\
                   & 4 & \tf 0.95870206828682745811 & \tf 0.20228209987458864481   \\\hline
\multirow{3}{*}{2} & 5 & \tf 1.17267316464601142810 & \tf 0.27222222222222222222   \\
                   & 6 & \tf 1.5                    & \tf 0.35555555555555555555   \\
                   & 7 & \tf 1.82732683535398857190 & \tf 0.27222222222222222222   \\\hline
\multirow{4}{*}{3} & 8 & \tf 2.04279465186386840500 & \tf 0.20648114717852759795   \\
                   & 9 & \tf 2.32101760363894084660 & \tf 0.34351885282147240205   \\
                   & 10& \tf 2.67898239636105915341 & \tf 0.34351885282147240205   \\
                   & 11& \tf 2.95720534813613159500 & \tf 0.20648114717852759795   \\\hline
\multirow{3}{*}{4} & 12& \tf 3.17267316464601142810 & \tf 0.27222222222222222222   \\
                   & 13& \tf 3.5                    & \tf 0.35555555555555555555   \\
                   & 14& \tf 3.82732683535398857190 & \tf 0.27222222222222222222   \\\hline
\multirow{4}{*}{5} & 15& \tf 4.04279465186386840500 & \tf 0.20648114717852759795   \\
                   & 16& \tf 4.32101760363894084660 & \tf 0.34351885282147240205   \\
                   & 17& \tf 4.67898239636105915341 & \tf 0.34351885282147240205   \\
                   & 18& \tf 4.95720534813613159500 & \tf 0.20648114717852759796   \\\hline
\end{tabular}
}
\end{center}
\end{table}

\begin{table}[!tb]
 \begin{center}
  \begin{minipage}{0.9\textwidth}
\caption{Nodes and weights for Gaussian quadrature rule (\ref{quadrature}) for $d=9$, $c=0$,
with uniform and normalized $N=7$ elements.
The nodes and weights are displayed with $20$ decimal digits, showing only the nodes and weights in the first element differ from
the asymptotic values.}\label{tab:90}
  \end{minipage}
\vspace{0.2cm}\\
\small{
\renewcommand{\arraystretch}{1.15}
\renewcommand{\tf}{\small}
\begin{tabular}{| c | c || l| l|}\hline
\multicolumn{2}{|c}{} & \multicolumn{2}{c|}{$d=9$, $c=0$, $N=7$, uniform, $\|\br\| = 7.53^{-24}$ } \\\hline
$\#$el. & $i$ &  \multicolumn{1}{c|}{$\tau_i$}  & \multicolumn{1}{c|}{$\omega_i$} \\\hline\hline
\multirow{5}{*}{1} & 1 & \tf 0.04769868312600247039 & \tf 0.12045551645375619589    \\
                   & 2 & \tf 0.23465028410411066421 & \tf 0.24335216912167039233    \\
                   & 3 & \tf 0.50845157638500186555 & \tf 0.28930851504236480624    \\
                   & 4 & \tf 0.78242184273917027918 & \tf 0.24378857317919028563    \\
                   & 5 & \tf 0.97223215910026017521 & \tf 0.13642855953635165325    \\\hline
\multirow{4}{*}{2} & 6 & \tf 1.11747233803526765358 & \tf 0.18923747814892349016    \\
                   & 7 & \tf 1.35738424175967745184 & \tf 0.27742918851774317651    \\
                   & 8 & \tf 1.64261575824032254816 & \tf 0.27742918851774317651    \\
                   & 9 & \tf 1.88252766196473234643 & \tf 0.18923747814892349016    \\\hline
\multirow{5}{*}{3} & 10& \tf 2.02843379420745745224 & \tf 0.13831195367835027217    \\
                   & 11& \tf 2.22138680297870775253 & \tf 0.24789494287337386576    \\
                   & 12& \tf 2.5                    & \tf 0.29425287356321839080    \\
                   & 13& \tf 2.77861319702129224747 & \tf 0.24789494287337386576    \\
                   & 14& \tf 2.97156620579254254776 & \tf 0.13831195367835027217    \\\hline
\multirow{4}{*}{4} & 15& \tf 3.11747233803526765358 & \tf 0.18923747814892349016    \\
                   & 16& \tf 3.35738424175967745184 & \tf 0.27742918851774317651    \\
                   & 17& \tf 3.64261575824032254816 & \tf 0.27742918851774317651    \\
                   & 18& \tf 3.88252766196473234643 & \tf 0.18923747814892349016    \\\hline
\end{tabular}
}
\end{center}
\end{table}

$C^0$ Nonics, $d=9$, $c=0$.
From Table~\ref{tab:90}, we again observe a very fast convergence to the asymptotic configuration ($N=\infty$) since
the $4.5$-node-per-element rule differs from the asymptotic counterpart only on the first element.

We emphasize here, that the asymptotic rules are exact only when the integration domain is a whole real line.
However, our rules which are exact and optimal over \emph{finite} domains quickly converge to their asymptotic counterparts
for low continuities.
That is, the nodes and weights on only first few boundary elements differ from the asymptotic values.
Therefore, one can use our nodes and weights computed on the first few boundary elements combined with the asymptotic
values for the intermediate intervals.

\subsection{Optimal quadrature rules for spaces with higher continuities}\label{ssec:ExNonUni}

\begin{figure}[!tb]
\vrule width0pt\hfill
 \begin{overpic}[width=.49\columnwidth,angle=0]{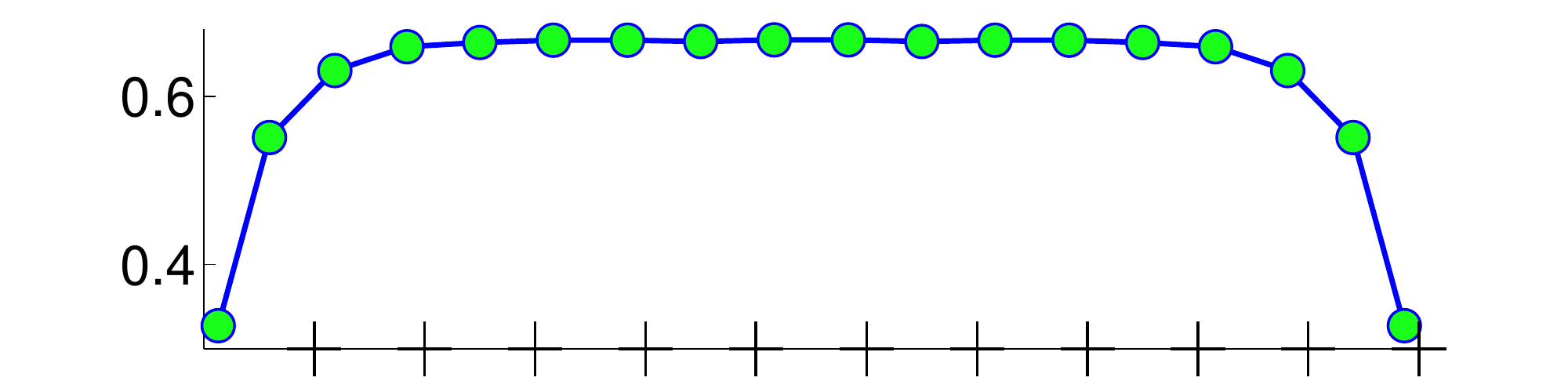}
    \put(90,30){\fcolorbox{gray}{white}{$d=5$, $c=2$}}
    \put(15,29){\fcolorbox{gray}{white}{$N=11$}}
    \put(82,-3){$x_{9}=11$}
	\end{overpic}
 \hfill
 \begin{overpic}[width=.49\columnwidth,angle=0]{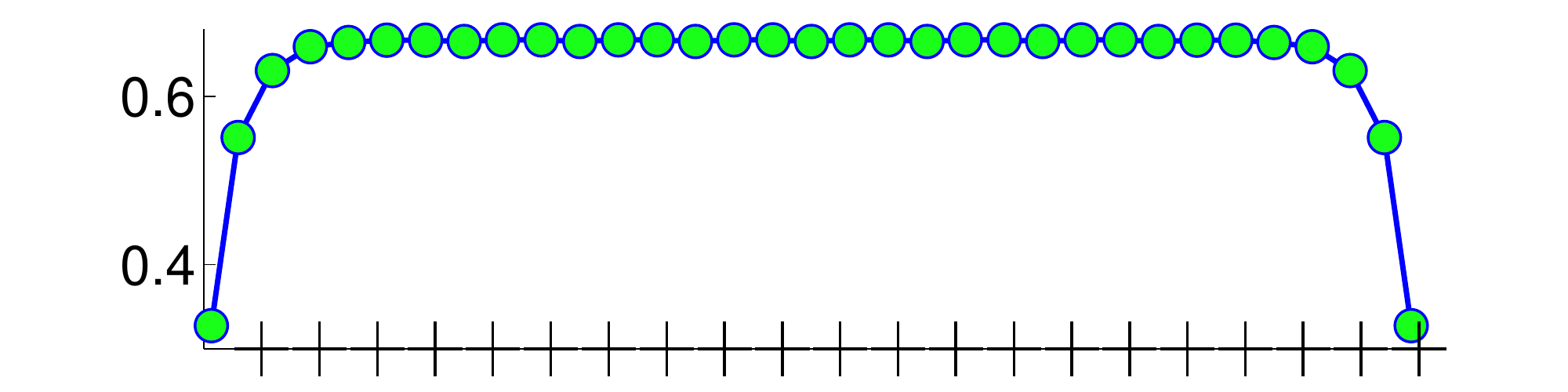}
    \put(75,29){\fcolorbox{gray}{white}{$N=21$}}
    \put(80,-3){$x_{21}=21$}
	\end{overpic}
 \hfill \vrule width0pt\\[1ex]
 \vrule width0pt\hfill
 \begin{overpic}[width=.89\columnwidth,angle=0]{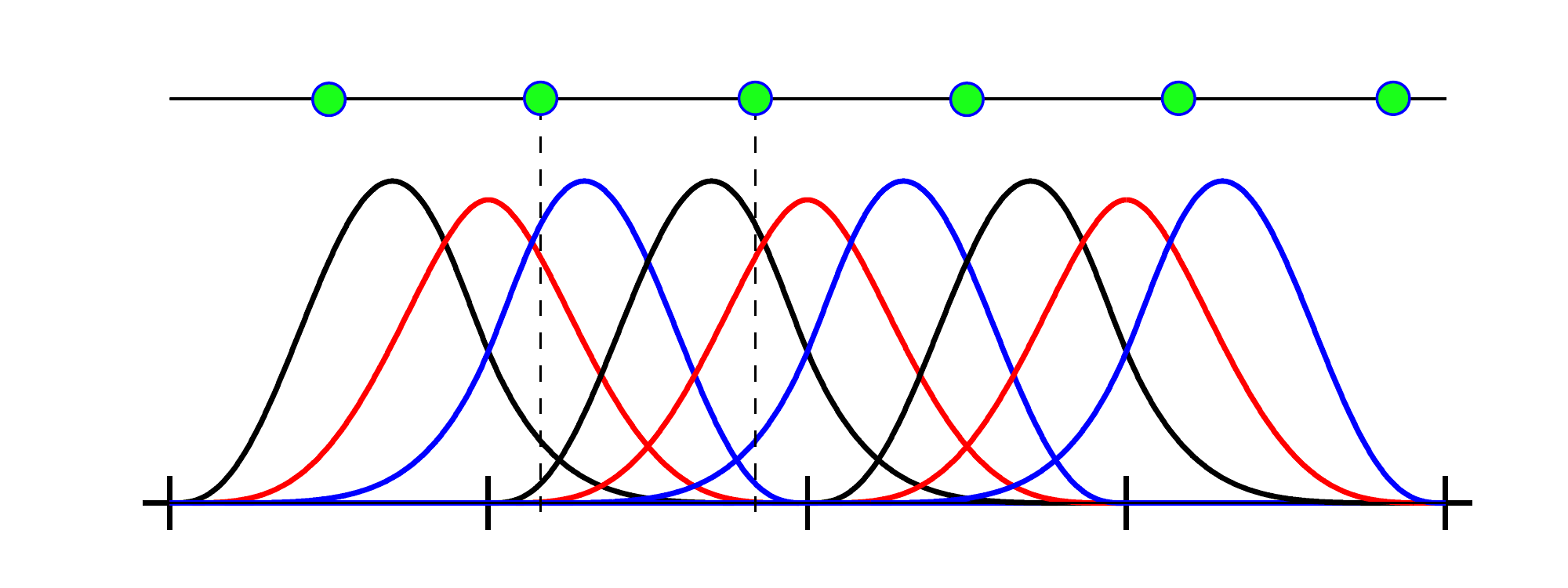}
    \put(8,1){$x_{i-1}$}
    \put(28,1){$x_i$}
    \put(48,-1){$d_1$}
    \put(31,-1){$d_1$}
    \put(51,1){$x_{i+1}$}
    \put(72,1){$x_{i+2}$}
    \put(30.5,4){$\scriptsize\underbrace{\makebox[0.4cm]{}}$}
    \put(47.5,4){$\scriptsize\underbrace{\makebox[0.4cm]{}}$}
    \put(4,30){\footnotesize $0.666$}
	\end{overpic}
 \hfill \vrule width0pt\\[0ex]
  \vrule width0pt\hfill
 \begin{overpic}[width=.93\columnwidth,angle=0]{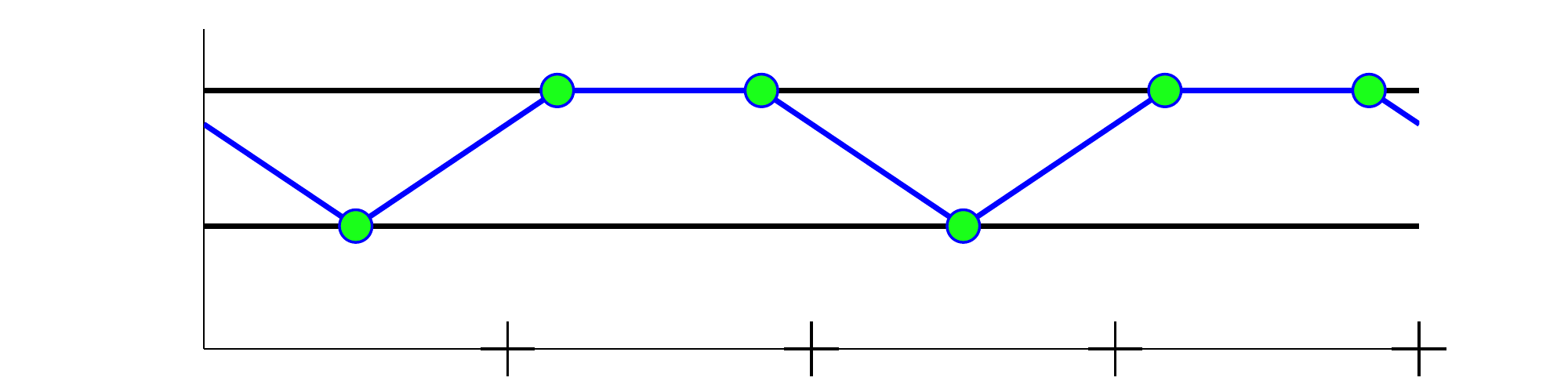}
     \put(2,10){\footnotesize $0.66554$}
     \put(2,19){\footnotesize $0.66723$}
	\end{overpic}
 \hfill \vrule width0pt\\[-2ex]
 \Acaption{1ex}{Top: Two optimal quadrature rules (green dots) for $C^2$-continuous quintic spline spaces for two different number of
 elements $N$ are shown.
  Middle: the asymptotic layout which requires $1.5$ nodes per element. Bottom: a zoom-in on the asymptotic weights.}\label{fig:d=5_k=2}
 \end{figure}

For quintic splines with $c=0$ and $c=1$ (see Tables~\ref{tab:50} and \ref{tab:51}),
the optimal nodes and weights differ from the asymptotic counterparts only on the first and on the first four boundary elements, respectively.
Higher continuity, however, drive the nodes and weights on significantly more boundary elements to differ from the asymptotic rules.

\begin{table}[!tb]
\begin{center}
\begin{minipage}{0.9\textwidth}
\caption{$1.5$-nodes-per-element Gaussian quadrature rule for $d=5$, $c=2$, with $N=31$ uniform elements over $[0,N]$.
Convergence to the asymptotic rule (\ref{eq:limit52}) is slow.
The nodes and weights on the first $15$ elements differ from asymptotic values by more than $16$ decimal digits.
}\label{tab:52}
  \end{minipage}
\vspace{0.2cm}\\
\small{
\renewcommand{\arraystretch}{1.15}
\renewcommand{\tf}{\small}
\begin{tabular}{| c | c || r| l|}\hline
\multicolumn{2}{|c}{} & \multicolumn{2}{c|}{$d=5$, $c=2$, $N=31$, uniform, $\|\br\| = 2.16^{-19}$ } \\\hline
$\#$el. & $i$ &  \multicolumn{1}{c|}{$\tau_i$}  & \multicolumn{1}{c|}{$\omega_i$} \\\hline\hline
\multirow{2}{*}{1} & 1 & \tf  0.13269173079619397468 & \tf 0.32754764055586624972    \\
                   & 2 & \tf  0.59618708045992614245 & \tf 0.55134370683162891514    \\\hline
\multirow{2}{*}{2} & 3 & \tf  1.18674528569337465654 & \tf 0.63102465144156062003    \\
                   & 4 & \tf  1.84114132687969848308 & \tf 0.65913151151996860606    \\\hline
                 3 & 5 & \tf  2.50087338091759294163 & \tf 0.66407995523024146993    \\\hline
\multirow{2}{*}{4} & 6 & \tf  3.16411265165773463104 & \tf 0.66694160636158860946    \\
                   & 7 & \tf  3.83609205795498629075 & \tf 0.66717534875798755120    \\\hline
                 5 & 8 & \tf  4.50000603299695540160 & \tf 0.66552621529342926507    \\\hline
\multirow{2}{*}{6} & 9 & \tf  5.16394446062823884836 & \tf 0.66722983983269530633    \\
                   & 10& \tf  5.83605695017807519279 & \tf 0.66723145226933460675    \\\hline
                 7 & 11& \tf  6.50000004155802779242 & \tf 0.66553624753029961958    \\\hline
\multirow{2}{*}{8} & 12& \tf  7.83605670832915656408 & \tf 0.66723182765303639270    \\
                   & 13& \tf  7.83605670832915656408 & \tf 0.66723183876013033758    \\\hline
                 9 & 14& \tf  8.50000000028626511144 & \tf 0.66553631663891425534    \\\hline
\multirow{2}{*}{10}& 15& \tf  9.16394329340371891384 & \tf 0.66723184134589561541    \\
                   & 16& \tf  9.83605670666322263574 & \tf 0.66723184142240485926    \\\hline
                11 & 17& \tf 10.50000000000197188624 & \tf 0.66553631711495679803    \\\hline
\multirow{2}{*}{12}& 18& \tf 11.16394329334871396842 & \tf 0.66723184144021644295    \\
                   & 19& \tf 11.83605670665174714652 & \tf 0.66723184144074346325    \\\hline
                13 & 20& \tf 12.50000000000001358298 & \tf 0.66553631711823593231    \\\hline
\multirow{2}{*}{14}& 21& \tf 13.16394329334833507662 & \tf 0.66723184144086615513    \\
                   & 22& \tf 13.83605670665166809954 & \tf 0.66723184144086978524    \\\hline
                15 & 23& \tf 14.50000000000000009292 & \tf 0.66553631711825851897    \\\hline
\end{tabular}
}
\end{center}
\end{table}

$C^2$ Quintics, $d=5$, $c=2$. 
Fig.~\ref{fig:d=5_k=2} shows two rules over finite domains for $N=11$ and $N=21$, and the asymptotic rule ($N=\infty$)
that requires $1.5$ nodes per element. In this scenario, one set of nodes becomes the midpoints (odd elements), while
the other nodes are posed pair-wise symmetrically with respect to elements' midpoints (even elements).
Due to the symmetry of the asymptotic basis functions, see Fig.~\ref{fig:d=5_k=2} middle, one can build a $(3\times 3)$ algebraic system.
The equations are again constraints of the exact integration, i.e. $\Q[D_i]=I[D_i]$, and the unknowns are $d_1$, $\omega_1$ and $\omega_2$.
Solving this asymptotic algebraic system, we obtain
\begin{equation}\label{eq:limit52}
\int_{\mathbb R}f(t) = \sum\limits_{i \in \mathbb Z} h( \omega_1 (f((2i+d_1)h) + f((2i+1-d_1)h)) + \omega_2 f(\frac{4i+3}{2}h))
\end{equation}
where $d_1$ is the smallest positive root of a univariate polynomial
\begin{equation}\label{eq:limit52roots}
-5+14x+150x^2-328x^3+164x^4
\end{equation}
and the weights $\omega_1$ and $\omega_2$ can be expressed in radicals, depending on $d_1$. Their numerical values are
\begin{equation}\label{eq:limit52dw}
\renewcommand{\arraystretch}{1.3}
\begin{array}{lcl}
d_1 & = & 0.83605670665166755138,  \\
w_1 & = & 0.66723184144087066164,  \\
w_2 & = & 0.66553631711825867672.
\end{array}
\end{equation}
The derived rule for $N=31$ is shown in Table~\ref{tab:52}.
For this $c=2$ case, we also observe convergence to the asymptotic rule (\ref{eq:limit52}),
however, the convergence is significantly slower as $15$ boundary elements are required to achieve a precision
of $16$ decimal digits.

$C^3$ Quintics, $d=5$, $c=3$. Even slower convergence is shown in Table~\ref{tab:53}.
The asymptotic rule requires one node per element (midpoint) and the limit weight is also one.
The rule is computed for the same number of elements as in the $c=2$ case, ($N=31$, cf. Table~\ref{tab:52}),
but here the nodes and weights on the $15$-th element meet the asymptotic values with only $10$ decimal digits.

\begin{table}[!tb]
\begin{center}
\begin{minipage}{0.9\textwidth}
\caption{One-node-per-element Gaussian quadrature rule for $d=5$, $c=3$, with $N=31$ uniform elements over $[0,N]$.
Observe a slow convergence to the asymptotic rule (elements' midpoints); the node on $15$-th element, $\tau_{16}$, meets the asymptotic value
only up to $10$ decimal digits.}\label{tab:53}
  \end{minipage}
\vspace{0.2cm}\\
\small{
\renewcommand{\arraystretch}{1.15}
\renewcommand{\tf}{\small}
\begin{tabular}{| c | c || r| l|}\hline
\multicolumn{2}{|c}{} & \multicolumn{2}{c|}{$d=5$, $c=3$, $N=31$, uniform, $\|\br\| = 1.29^{-17}$ } \\\hline
$\#$el. & $i$ &  \multicolumn{1}{c|}{$\tau_i$}  & \multicolumn{1}{c|}{$\omega_i$} \\\hline\hline
\multirow{2}{*}{1} & 1 & \tf  0.15502116527535364186 & \tf 0.38432797072681461548    \\
                   & 2 & \tf  0.72629241316800832428 & \tf 0.72643787278541483934    \\\hline
                 2 & 3 & \tf  1.55795146959426904181 & \tf 0.91383651913294099189    \\\hline
                 3 & 4 & \tf  2.51229503417716167402 & \tf 0.98033035056323515086    \\\hline
                 4 & 5 & \tf  3.50243819000163217876 & \tf 0.99603266582760659143    \\\hline
                 5 & 6 & \tf  4.50047612975551679586 & \tf 0.99922264148268942270    \\\hline
                 6 & 7 & \tf  5.50009269292182007059 & \tf 0.99984856434714317220    \\\hline
                 7 & 8 & \tf  6.50001803460249141356 & \tf 0.99997053247313114324    \\\hline
                 8 & 9 & \tf  7.50000350845339606943 & \tf 0.99999426724168815795    \\\hline
                 9 & 10& \tf  8.50000068251932617787 & \tf 0.99999888476863222037    \\\hline
                 10& 11& \tf  9.50000013277376488047 & \tf 0.99999978304847734739    \\\hline
                 11& 12& \tf 10.50000002582909436870 & \tf 0.99999995779540275690    \\\hline
                 12& 13& \tf 11.50000000502464245581 & \tf 0.99999999178972924635    \\\hline
                 13& 14& \tf 12.50000000097741787292 & \tf 0.99999999840273582518    \\\hline
                 14& 15& \tf 13.50000000018988049624 & \tf 0.99999999968884882908    \\\hline
                 15& 16& \tf 14.50000000003559360097 & \tf 0.99999999993726772687    \\\hline
\end{tabular}
}
\end{center}
\end{table}

\subsection{Optimal quadrature rules for spaces with non-uniform knots and continuities}\label{ssec:ExNonUni}

\begin{figure}[!tb]
\vrule width0pt\hfill
 \begin{overpic}[width=.82\columnwidth,angle=0]{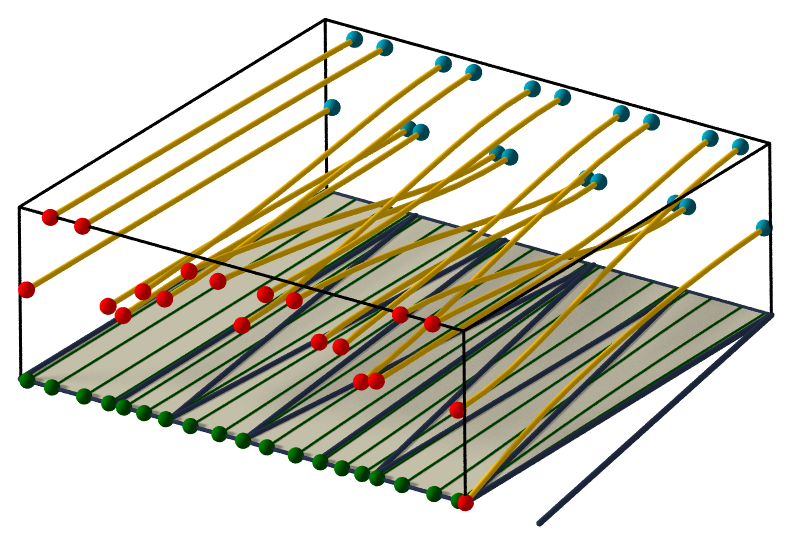}
   \put(2,10){$(\bx,\bm)$}
   \put(70,60){$(\xxt,\bmt)$}
   \put(-6,44){\fcolorbox{gray}{white}{\small $t=1$}}
   \put(35,67){\fcolorbox{gray}{white}{\small $t=0$}}
   \put(2,5){\small $\bx=(0,\frac{5}{24},\frac{1}{3},\frac{1}{2},\frac{2}{3},\frac{19}{24},1,\frac{7}{6})$}
   \put(2,0){\small $\bm=(8,7,5,6,5,7,8,2)$}
   \put(55,70){\small $\xxt=(0,\frac{1}{5},\frac{2}{5},\frac{3}{5},\frac{4}{5},1)$}
   \put(55,65){\small $\bmt=(8,8,8,8,8,8)$}
\end{overpic}
\hfill \vrule width0pt\\[2ex]
\vrule width0pt\hfill
 \begin{overpic}[width=.85\columnwidth,angle=0]{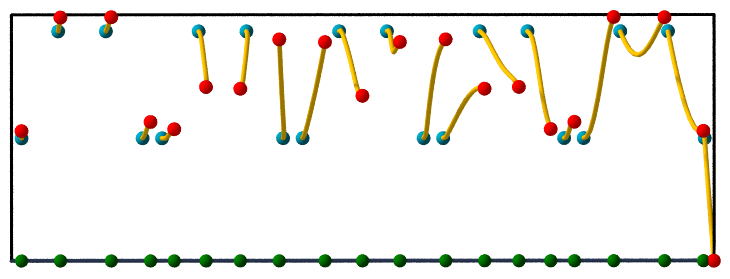}
 \put(-2,0){\small $0$}
 \put(100,0){\small $1$}
 \put(95,38.5){\small $[1,0.07]$}
\end{overpic}\hfill \vrule width0pt\\[-3ex]
\Acaption{1ex}{The evolution of Gaussian quadrature rule (\ref{quadrature}) for $d=7$, $N=6$ ($n=5$) is shown.
Top: starting from the source space $\St^{n,d}_{\xxt,\bmt}$ generated above $(\xxt,\bmt)$,
the source rule (blue) is evolved into the target rule (red) via a geodesic knot transformation (dark blue).
The target space $S^{N,d}_{\bx,\bm}$ is defined above a knot vector with non-uniform partition $\bx$
and non-uniform multiplicities $\bm$.
 Bottom: Front view of the evolution with the target rule ($t=1$) highlighted in red.
 The dimensions of spaces $\St^{5,7}_{\xxt,\bmt}$ and $S^{6,7}_{\bx,\bm}$ on $[a,b]=[0,1]$ are $40$ and $38$, respectively.
 Therefore two knots (one double knot) are pulled outside $[a,b]$, which forces $\tau_{20} = b$ and $\omega_{20} = 0$ at the target time $t=1$.
}\label{fig:VariousTargetsN6}
 \end{figure}

So far, we focused only on the target spaces with uniform knots, however, the homotopy continuation
method is well-suited for non-uniform knot distributions.
The evolution of the Gaussian rule for a particular target spline space with both non-uniform knot partitions
and multiplicities is shown in Fig.~\ref{fig:VariousTargetsN6}.
The continuity varies across the elements (from $ C^0$ to $C^2$).
The corresponding nodes and weights derived by our algorithm are listed in Table~\ref{tab:nonuni},
showing also the errors of the rule measured by (\ref{eq:Error}).
%

\begin{table}[!tb]
\begin{center}
\begin{minipage}{0.9\textwidth}
\caption{Gaussian quadrature rule for $d=7$ with varying continuities and $N=6$ non-uniform elements over $[0,1]$ shown in Fig.~\ref{fig:VariousTargetsN6}.
Due to symmetry, only the first ten values on $[0,0.5]$ are shown.
}\label{tab:nonuni}
  \end{minipage}
\vspace{0.2cm}\\
\small{
\renewcommand{\arraystretch}{1.15}
\renewcommand{\tf}{\small}
\begin{tabular}{| c | c || l| l|}\hline
\multicolumn{2}{|c}{} & \multicolumn{2}{c|}{$d=7$, $c$ varies, $N=6$, non-uniform, $\|\br\| = 3.61^{-17}$ } \\\hline
$\#$el. & $i$ &  \multicolumn{1}{c|}{$\tau_i$}  & \multicolumn{1}{c|}{$\omega_i$} \\\hline\hline
\multirow{4}{*}{1} & 1 & \tf 0.01475556054370093982 & \tf 0.03696325520404057241  \\
                   & 2 & \tf 0.07013752839667597547 & \tf 0.06930847014854596998 \\
                   & 3 & \tf 0.14242843659730234676 & \tf 0.06938624349804886061  \\
                   & 4 & \tf 0.19834373082551928914 & \tf 0.03955024610641599980  \\\hline
\multirow{3}{*}{2} & 5 & \tf 0.23209282669863971447 & \tf 0.03749098184955721430  \\
                   & 6 & \tf 0.27734376491429861094 & \tf 0.04948821979402789722  \\
                   & 7 & \tf 0.32617242559157607738 & \tf 0.04899222891733174032  \\\hline
\multirow{3}{*}{3} & 8 & \tf 0.38159151883896678342 & \tf 0.06305153787491328839  \\
                   & 9 & \tf 0.44667150881630100465 & \tf 0.06228144839443712894  \\
                   & 10& \tf 0.50000000000000010463 & \tf 0.04697473642536305323  \\\hline
\end{tabular}
}
\end{center}
\end{table}


 \begin{table}[!tb]
 \begin{center}
  \begin{minipage}{0.9\textwidth}
\caption{Overview of the derived uniform Gaussian rules for various $d$ and $c$ and their correlation
with the asymptotic counterparts.
The integers in every box refer to the boundary elements where the rules differ from their asymptotic counterparts
by more than double machine precision. Where applicable, tables that contain the boundary nodes and weights
as well as the corresponding asymptotic rules are shown.  }\label{tab:Overview}
  \end{minipage}
\vspace{0.2cm}\\
\small{
\renewcommand{\arraystretch}{1.15}
\renewcommand{\tf}{\small}
\begin{tabular}{| c || c | c | c | c |}\hline 
\rotatebox{0}{}
 \makebox[0.6cm]{} $c$ &  $0$  & $1$ &  $2$  & $3$\\
 $d$  &  &  &  & \\\hline\hline
 \multirow{2}{*}{$5$} & 1 & 4 & 15 & $\gg$15 \\
    & Tab.~\ref{tab:50}, (\ref{eq:limit50}) &  Tab.~\ref{tab:51}, (\ref{eq:InfyWeights51}) & Tab.~\ref{tab:52}, (\ref{eq:limit52}) & Tab.~\ref{tab:53} \\\hline
 \multirow{2}{*}{$7$} & 1 & 4 & 14 & $\gg$14 \\
    & Tab.~\ref{tab:70} & Tab.~\ref{tab:71}, (\ref{eq:limit71}) & Appendix & Appendix \\\hline
 \multirow{2}{*}{$9$} & 1 & 4  & &  \\
    & Tab.~\ref{tab:90}  & Tab.~\ref{tab:91}, (\ref{eq:limit91}) & &  \\\hline
\end{tabular}
}
\end{center}
\end{table}

\subsection{Summary of the derived rules}\label{ssec:Summary}

Table~\ref{tab:Overview} shows the summary of the derived optimal rules for spline spaces of various degrees
and continuities over finite domains. Where applicable, it also shows the corresponding asymptotic formulae and
displays the number of boundary elements which differ, up to machine precision, from the asymptotic values.
From Table~\ref{tab:Overview}, we conjecture that the convergence to the asymptotic rules
depends on continuity rather than on polynomial degree. While for low continuities ($C^0$ and $C^1$) only a few boundary elements differ from the asymptotic
values, for higher continuities one needs significantly more elements and it is therefore less convenient
to be combined with the asymptotic counterparts.

The results in all the tables show the errors of the derived rules using the Euclidean metric on the vector of residues, $\|\br\|$,
when applied to the basis functions. Even though this already measures quality of the rule, it is informative to see
the actual error of the rule when applied to a function from the
space under consideration. We tested our rules against random functions from the particular spline spaces.
Few examples are shown in the Appendix.
In all our tests, the exact integrals and our quadrature rules match up to sixteenth decimal
digit.

\section{Conclusion}\label{sec:conl}

We introduce a new methodology to compute Gaussian quadrature rules for
spline spaces that are frequently used in Galerkin discretizations
when building mass and stiffness matrices using isogeometric analysis.
The rules use the minimum
number of quadrature points, and are exact when integrating over finite domains.
Starting with the classical Gaussian quadrature for polynomials on a specifically designed mesh,
the optimal rule is interpreted as a zero of a particular polynomial system.
We use homotopy continuation to
numerically trace the Gaussian quadrature rule for a target space as the knot vector,
and consequently the whole spline space, is continuously modified.
We show that the methodology handles scenarios where the source space and the target
space possess different numbers of optimal quadrature points.

The numerical examples demonstrate that the methodology is able to find optimal quadrature rules on spline spaces of various degrees
and continuities. We show numerically that our rules satisfy very tight numerical thresholds and
therefore are exact to machine precision.
We discuss the asymptotic behavior as $N$ grows
and show the convergence to the asymptotic rules (for some spaces derived in \cite{Hughes-2010}),
where the integration domain is a whole real line.
For lower continuity cases, see Table~\ref{tab:Overview}, we numerically show a fast convergence to the asymptotic layout
and therefore only the nodes and weights on the first few elements differ from the repetitive asymptotic pattern.
Consequently, a combination of the derived nodes and weights together with the proper asymptotic rule forms a quasi-repetitive pattern that
is convenient and easy to implement. We believe the optimal rules derived here
could become standard in IGA-oriented software such as PetIGA \cite{petiga},
which so far uses classical polynomial Gaussian quadratures that are suboptimal.

We also emphasize that besides the minimum number of quadrature points,
all the rules have positive weights. This fact makes them preferable choice from the point of view of numerical stability,
e.g., when compared to \cite{Hughes-2012} where additional quadrature points are added, some of them having negative weights.


As future work, we plan to investigate
Gauss-Lobatto rules for particular spline spaces since the homotopy-based approach seems to be a very elegant tool.
Likewise, optimal rules for spline spaces of even degree remain a very challenging problem.

\section*{References}\label{sec:ref}

\bibliographystyle{plain}
\bibliography{StiffnessMatrixArXiv}

\clearpage

\section*{Appendix}

Gaussian quadrature rules for $C^2$ and $C^3$ septic spline spaces are shown in Tables \ref{tab:72} and \ref{tab:73},
respectively. We observe a slow convergence to the asymptotic counterparts. With fixed number of elements $N=31$,
for $c=2$, the nodes and weights on the first fourteen boundary elements differ from the limit values by more than
machine precision in Table~\ref{tab:72}. On the contrary, for $c=3$ in Table~\ref{tab:73},
the nodes and weights on the $14$-th element match the asymptotic values by only $12$ decimal digits.

\begin{figure}[!b]
 \vrule width0pt\hfill
 \begin{overpic}[width=.49\columnwidth,angle=0]{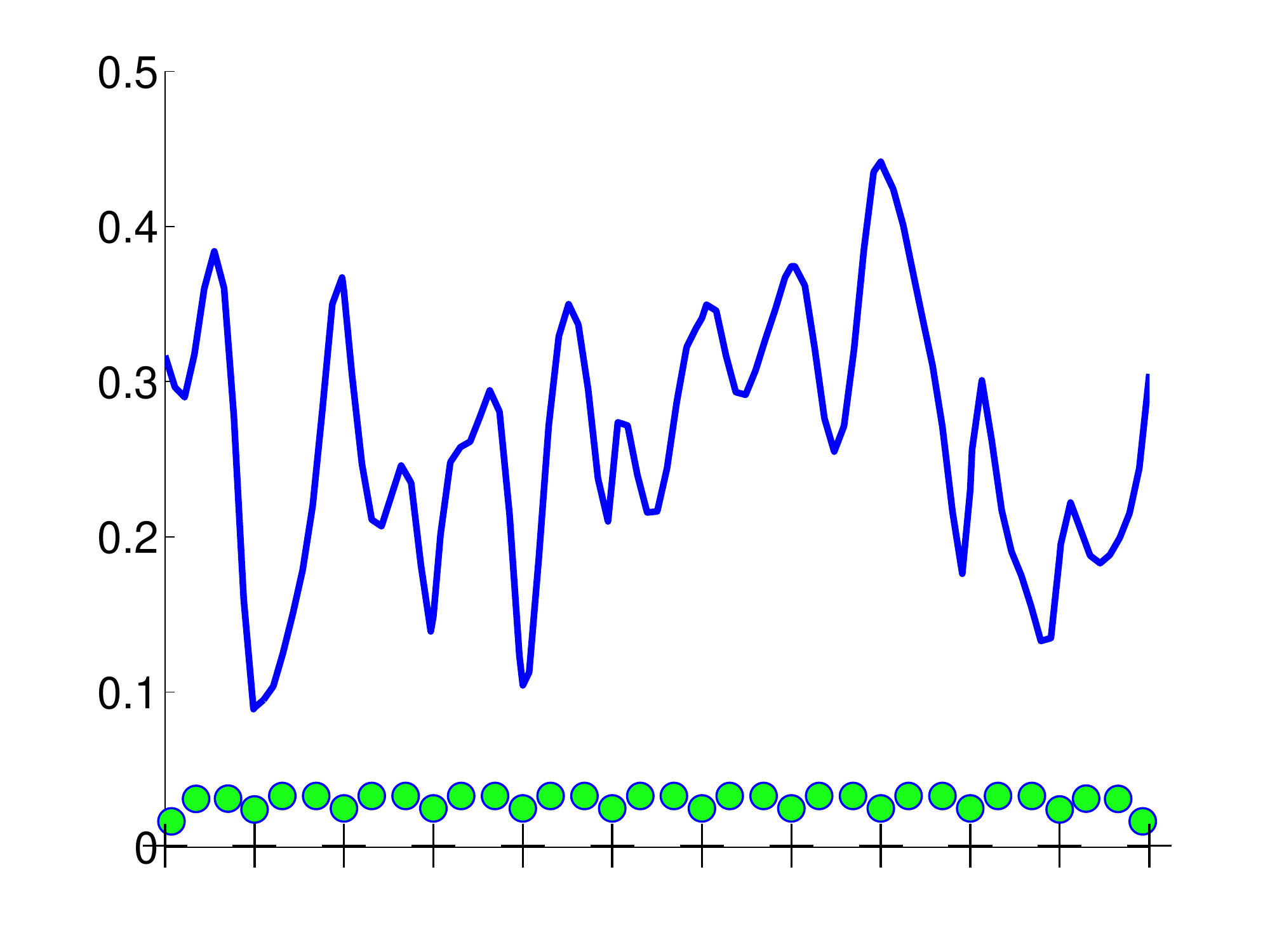}
    \put(15,65){\fcolorbox{gray}{white}{$\varepsilon = 2.3^{-17}$}}
    \put(55,80){\fcolorbox{gray}{white}{$d=7$, $c=1$, $N=11$, $\card(\XX_N) = 76$, $r=4$}}
    \put(15,3){$x_1$}
    \put(82,3){$x_{11}=1$}
	\end{overpic}
 \hfill
 \begin{overpic}[width=.49\columnwidth,angle=0]{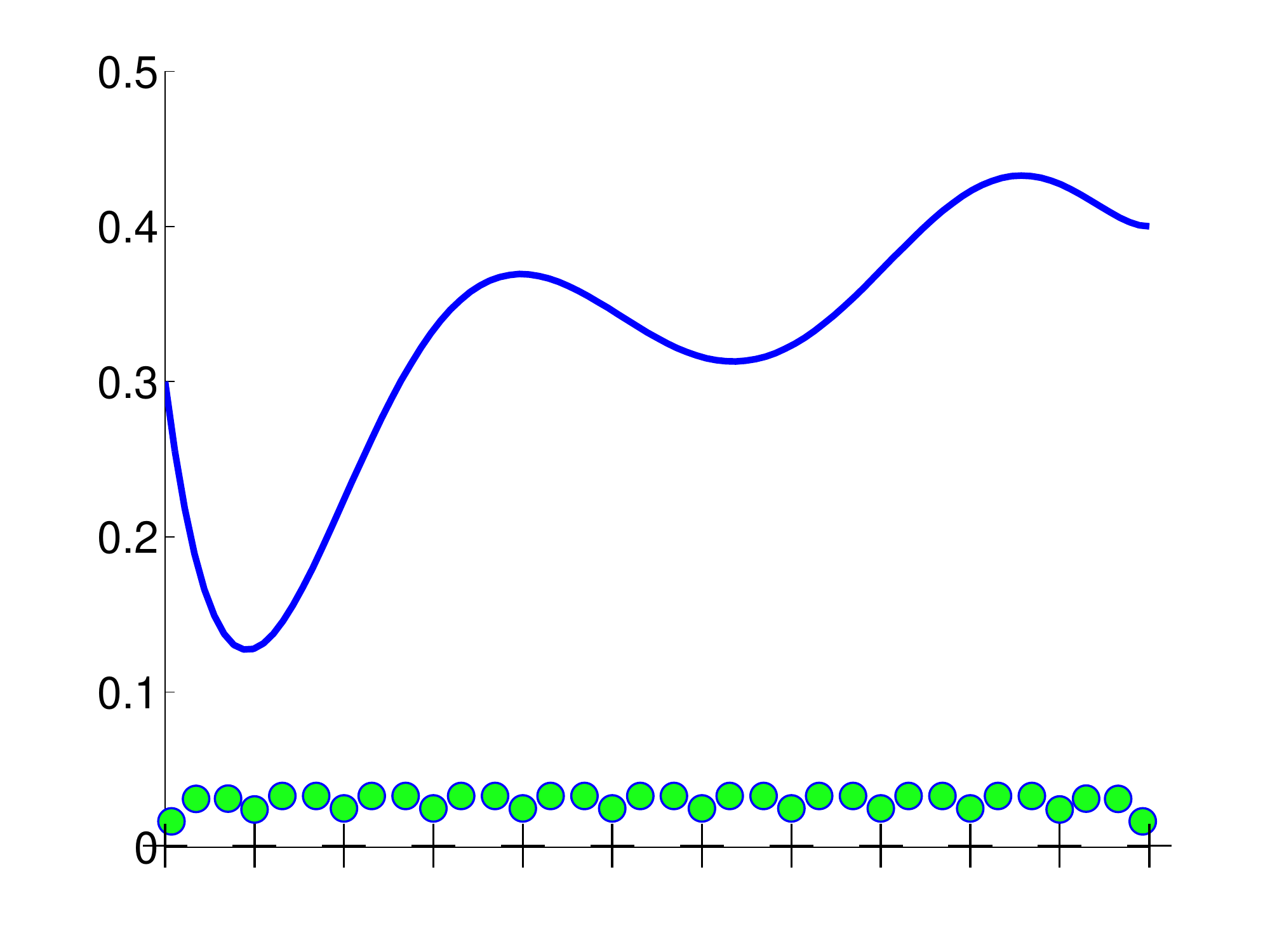}
    \put(15,65){\fcolorbox{gray}{white}{$\varepsilon = 5.7^{-17}$}}
    \put(15,3){$x_1$}
    \put(80,3){$x_{11}=1$}
    \put(25,20){\small $\bb = (\frac{3}{10},-\frac{2}{5}, 1, \frac{3}{5},\frac{1}{2},\frac{4}{5},\frac{2}{5},\frac{2}{5})$}
	\end{overpic}
 \hfill \vrule width0pt\\[-5ex]
 \Acaption{1ex}{Numerical tests for uniform $d=7$, $c=1$, $N=11$. Two functions (blue) from the spline space
 over $N=11$ elements are shown together with the integration errors $\varepsilon$ between our rule (green dots) and exact integration.
 A function with random spline basis coefficients (left)
 and a degree seven polynomial with a vector $\bb$ of its B\'ezier coefficients (right) are displayed.}\label{fig:NumTest71}
 \end{figure}

\begin{figure}[!b]
 \vrule width0pt\hfill
 \begin{overpic}[width=.49\columnwidth,angle=0]{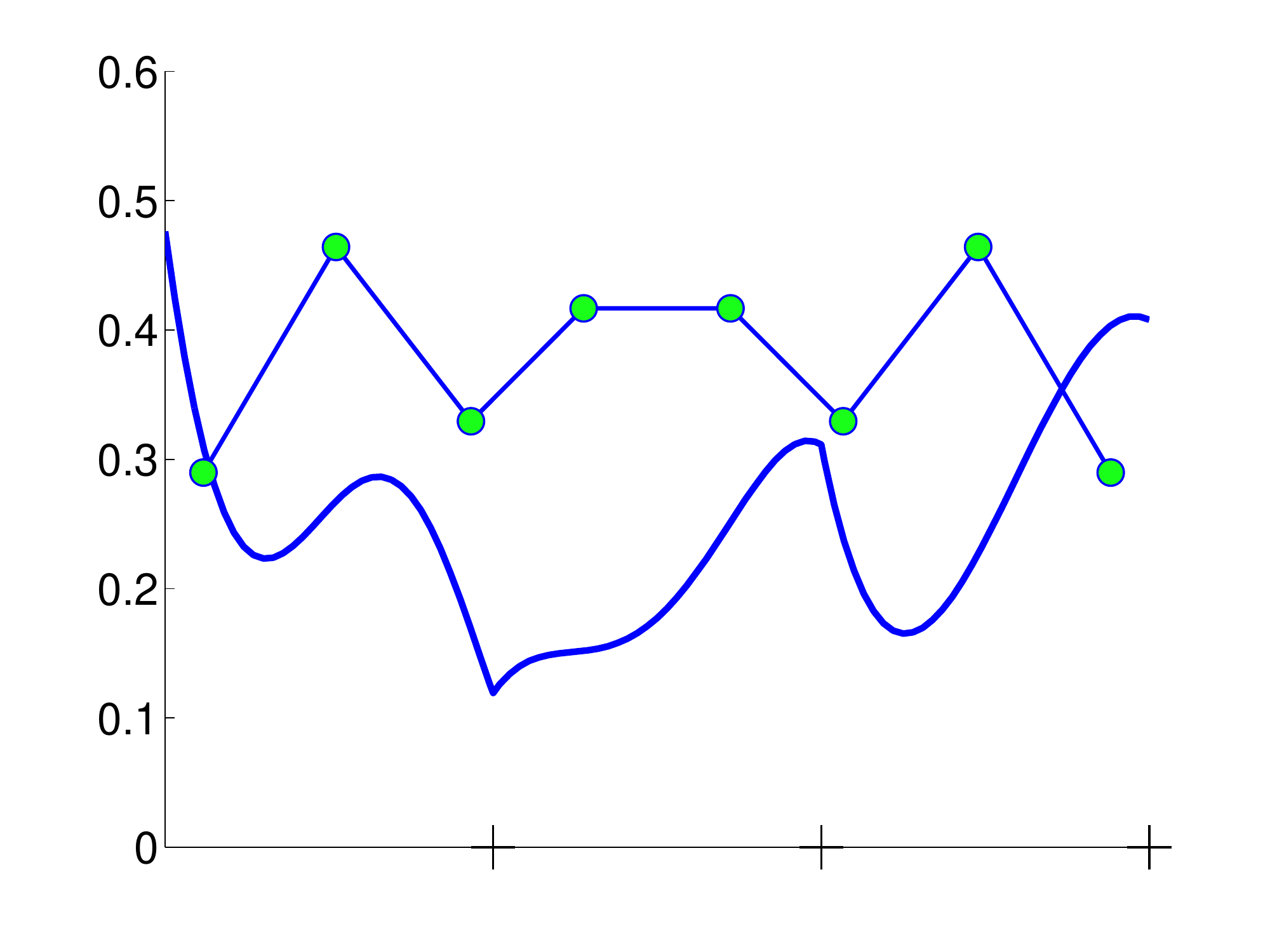}
    \put(15,65){\fcolorbox{gray}{white}{$N=3$, $r=2$}}
    \put(60,15){\fcolorbox{gray}{white}{$\varepsilon = 4.1^{-18}$}}
    \put(82,3){$x_{3}=3$}
	\end{overpic}
 \hfill
 \begin{overpic}[width=.49\columnwidth,angle=0]{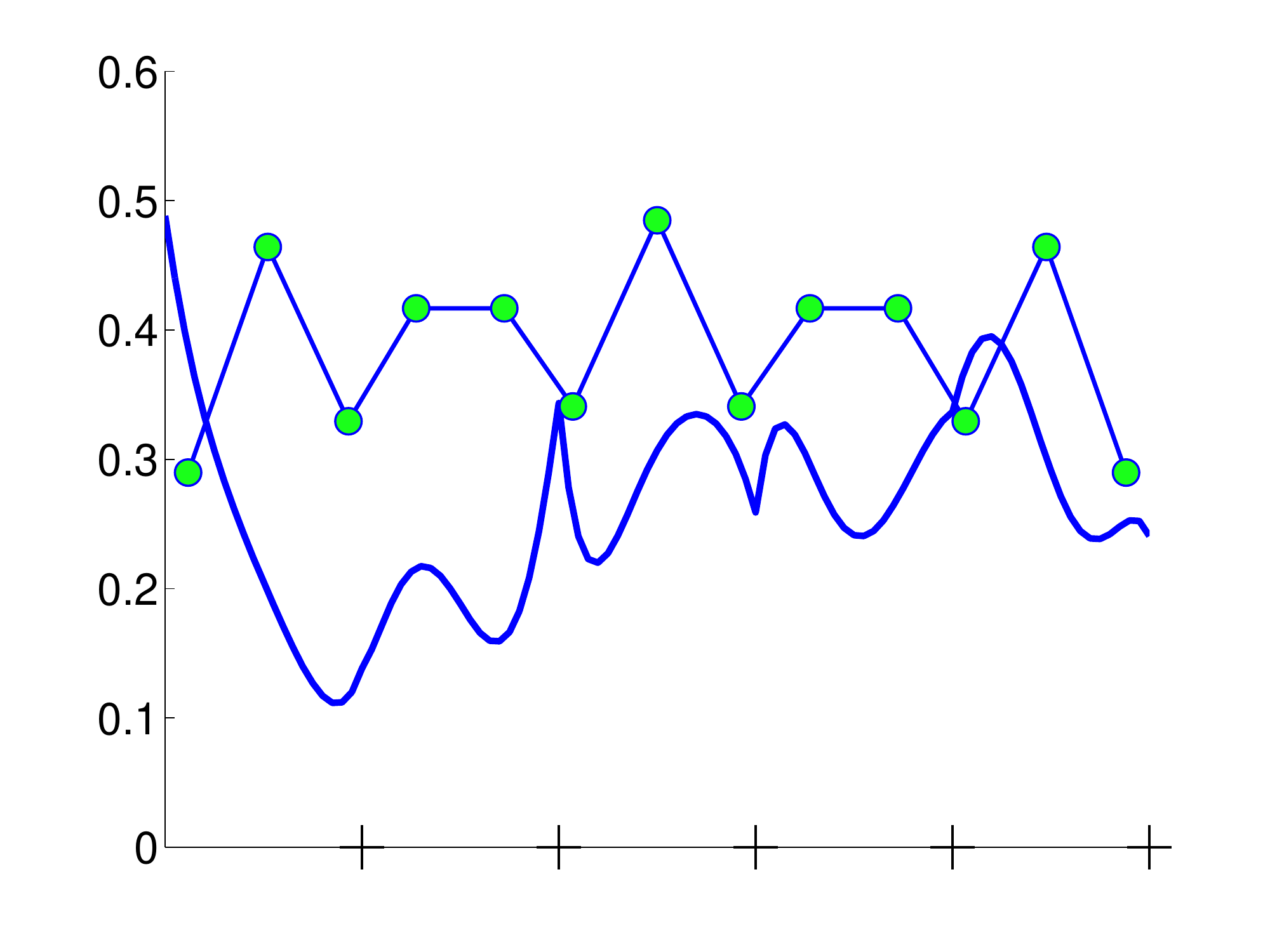}
    \put(60,65){\fcolorbox{gray}{white}{$N=5$, $r=2$}}
    \put(60,15){\fcolorbox{gray}{white}{$\varepsilon = 2.7^{-17}$}}
    \put(-15,70){\fcolorbox{gray}{white}{$d=5$, $c=0$}}
    \put(80,3){$x_{5}=5$}
	\end{overpic}
 \hfill \vrule width0pt\\[-5ex]
 \Acaption{1ex}{Optimal quadrature rules for $C^0$-continuous quintic spline spaces with $N$ uniform elements.
 Our quadrature rules are visualized by the green dots and the errors between our rule and
 exact integration are shown. The functions (blue) were generated by random sampling of the spline coefficients
 from the uniform distribution on $[0,0.5]$.}\label{fig:NumTest50}
 \end{figure}

Figure~\ref{fig:NumTest71} shows results of a numerical test of our quadrature rule
when applied to two functions from the spaces under consideration for $d=7$, $c=1$, $N=11$.
The first function was generated randomly by sampling its spline coefficients,
the second one is a polynomial of degree seven.
Figure~\ref{fig:NumTest50} shows two random function for $d=5$, $c=0$, and $N=3,5$.
The actual error $\varepsilon$ between the derived rules and exact (symbolic) integrals are shown.

\begin{table}[!h]
 \begin{center}
  \begin{minipage}{0.9\textwidth}
\caption{Nodes and weights for Gaussian quadrature rule for $C^2$ septics, $d=7$, $c=2$, with uniform and normalized $N=31$ elements.
The nodes and weights are displayed with $20$ decimal digits, showing the nodes and weights in the first $14$ elements which differ from
the asymptotic values by more than $16$ decimal digits.}\label{tab:72}
  \end{minipage}
\vspace{0.2cm}\\
\small{
\renewcommand{\arraystretch}{1.15}
\renewcommand{\tf}{\small}
\begin{tabular}{| c | c || r| l|}\hline
\multicolumn{2}{|c}{} & \multicolumn{2}{c|}{$d=7$, $c=2$, $N=31$, uniform, $\|\br\| = 2.92^{-19}$ } \\\hline
$\#$el. & $i$ &  \multicolumn{1}{c|}{$\tau_i$}  & \multicolumn{1}{c|}{$\omega_i$} \\\hline\hline
\multirow{3}{*}{1} & 1 & \tf 0.07602281159668294261 & \tf 0.19047916188601100653    \\
                   & 2 & \tf 0.36176994290688769251 & \tf 0.35825473257063950190    \\
                   & 3 & \tf 0.73905374789752879161 & \tf 0.37182997478207990632    \\\hline
\multirow{3}{*}{2} & 5 & \tf 1.09855791031238813896 & \tf 0.36855197226751600514    \\
                   & 6 & \tf 1.50262223080776729934 & \tf 0.42715636781540455283   \\
                   & 7 & \tf 1.90950254570512762462 & \tf 0.37775194547562865077    \\\hline
\multirow{2}{*}{3} & 8 & \tf 2.29132365358414524361 & \tf 0.40607590839141959194     \\
                   & 9 & \tf 2.70890362556701209765 & \tf 0.40619403944394371367    \\\hline
\multirow{3}{*}{4} & 12& \tf 3.09130818926984030943 & \tf 0.37899102728460633495    \\
                   & 13& \tf 3.50000479301042700738 & \tf 0.42943959528926707300    \\
                   & 14& \tf 3.90870701779525770041 & \tf 0.37900947826151803192   \\\hline
\multirow{2}{*}{5} & 15& \tf 4.29115702707950753516 & \tf 0.40626597892248379641    \\
                   & 16& \tf 4.70884338733068500007 & \tf 0.40626619496240698197   \\\hline
\multirow{3}{*}{6} & 17& \tf 5.09129446302025060221 & \tf 0.37901174209623557803    \\
                   & 18& \tf 5.50000000873421055658 & \tf 0.42944377906050160659   \\
                   & 19& \tf 5.90870556469302714690 & \tf 0.37901177572516939951    \\\hline
\multirow{2}{*}{7} & 20& \tf 6.29115672321619455811 & \tf 0.40626632609753336879   \\
                   & 21& \tf 6.70884327753897390813 & \tf 0.40626632649121921345    \\\hline
\multirow{3}{*}{8} & 22& \tf 7.09129443800543126056 & \tf 0.37901177985050887499   \\
                   & 23& \tf 7.50000000001591609819 & \tf 0.42944378668453946421    \\
                   & 24& \tf 7.90870556204506984624 & \tf 0.37901177991178992443   \\\hline
\multirow{2}{*}{9} & 25& \tf 8.29115672266247252821 & \tf 0.40626632673018308189   \\
                   & 26& \tf 8.70884327733890359334 & \tf 0.40626632673090048417   \\\hline
\multirow{3}{*}{10}& 27& \tf 9.09129443795984747612 & \tf 0.37901177991930741072   \\
                   & 28& \tf 9.50000000000002900345 & \tf 0.42944378669843252704   \\
                   & 29& \tf 9.90870556204024455058 & \tf 0.37901177991941908139   \\\hline
\multirow{2}{*}{11}& 30& \tf10.29115672266146349664 & \tf 0.40626632673133594109   \\
                   & 31& \tf10.70884327733853901102 & \tf 0.40626632673133724839   \\\hline
\multirow{3}{*}{12}& 32& \tf11.09129443795976441012 & \tf 0.37901177991943278029   \\
                   & 33& \tf11.50000000000000005285 & \tf 0.42944378669845784397    \\
                   & 34& \tf11.90870556204023575758 & \tf 0.37901177991943298378   \\\hline
\multirow{2}{*}{13}& 35& \tf12.29115672266146165791 & \tf 0.40626632673133804191   \\
                   & 36& \tf12.70884327733853834666 & \tf 0.40626632673133804429   \\\hline
\multirow{3}{*}{14}& 37& \tf13.09129443795976425875 & \tf 0.37901177991943300875   \\
                   & 38& \tf13.50000000000000000010 & \tf 0.42944378669845789010   \\
                   & 39& \tf13.90870556204023574156 & \tf 0.37901177991943300912   \\\hline
\end{tabular}
}
\end{center}
\end{table}

\begin{table}[!tb]
 \begin{center}
  \begin{minipage}{0.9\textwidth}
\caption{Two-nodes-per-element Gaussian quadrature rule
for $d=7$, $c=3$, with $N=31$ uniform elements over $[0,N]$.
The nodes and weights are displayed with $20$ decimal digits, showing a slow convergence to the asymptotic rule.
The nodes and weights on the $14$-th element meet the asymptotic values by only $12$ decimal digits.}\label{tab:73}
  \end{minipage}
\vspace{0.2cm}\\
\small{
\renewcommand{\arraystretch}{1.15}
\renewcommand{\tf}{\small}
\begin{tabular}{| c | c || l| l|}\hline
\multicolumn{2}{|c}{} & \multicolumn{2}{c|}{$d=7$, $c=3$, $N=31$, uniform, $\|\br\| = 6.41^{-18}$ } \\\hline
$\#$el. & $i$ &  \multicolumn{1}{c|}{$\tau_i$}  & \multicolumn{1}{c|}{$\omega_i$} \\\hline\hline
\multirow{2}{*}{1} & 1 & \tf 0.08173435864484871477 & \tf 0.20488801504260514819   \\
                   & 2 & \tf 0.39001985639208382229 & \tf 0.38827700271568336355   \\
                   & 3 & \tf 0.80983038609710312021 & \tf 0.43593835564819134630   \\\hline
\multirow{2}{*}{2} & 4 & \tf 1.26440542030435148338 & \tf 0.48028254678989050587   \\
                   & 5 & \tf 1.75806373790131509116 & \tf 0.49335481875705423628   \\\hline
\multirow{2}{*}{3} & 6 & \tf 2.24908303781589110505 & \tf 0.49809212627244622694   \\
                   & 7 & \tf 2.75295470311221074603 & \tf 0.49940844564325211137   \\\hline
\multirow{2}{*}{4} & 8 & \tf 3.24767613649844940900 & \tf 0.49983160560169774755   \\
                   & 9 & \tf 3.75250264266952998560 & \tf 0.49994819564193698306   \\\hline
\multirow{2}{*}{5} & 10& \tf 4.24755256119432444550 & \tf 0.49998526412553038651   \\
                   & 11& \tf 4.75246307519981025547 & \tf 0.49999546981721281010   \\\hline
\multirow{2}{*}{6} & 12& \tf 5.24754175195116296507 & \tf 0.49999871146226766302   \\
                   & 13& \tf 5.75245961526710334230 & \tf 0.49999960389454759273   \\\hline
\multirow{2}{*}{7} & 14& \tf 6.24754080680196649870 & \tf 0.49999988733480559009   \\
                   & 15& \tf 6.75245931274229531708 & \tf 0.49999996536611859321   \\\hline
\multirow{2}{*}{8} & 16& \tf 7.24754072416172307615 & \tf 0.49999999014900950818   \\
                   & 17& \tf 7.75245928629074166596 & \tf 0.49999999697175439787   \\\hline
\multirow{2}{*}{9} & 18& \tf 8.24754071693599533383 & \tf 0.49999999913866953032   \\
                   & 19& \tf 8.75245928397792556542 & \tf 0.49999999973522255412   \\\hline
\multirow{2}{*}{10}& 20& \tf 9.24754071630420715775 & \tf 0.49999999992468877663   \\
                   & 21& \tf 9.75245928377570235808 & \tf 0.49999999997684894010   \\\hline
\multirow{2}{*}{11}& 22& \tf10.24754071624896617647 & \tf 0.49999999999341509381   \\
                   & 23& \tf10.75245928375802078566 & \tf 0.49999999999797576575   \\\hline
\multirow{2}{*}{12}& 24& \tf11.24754071624413613038 & \tf 0.49999999999942424266   \\
                   & 25& \tf11.75245928375647478108 & \tf 0.49999999999982300921   \\\hline
\multirow{2}{*}{13}& 26& \tf12.24754071624371381085 & \tf 0.49999999999994965813   \\
                   & 27& \tf12.75245928375633960475 & \tf 0.49999999999998452465   \\\hline
\multirow{2}{*}{14}& 28& \tf13.24754071624367688495 & \tf 0.49999999999999559831   \\
                   & 29& \tf13.75245928375632778547 & \tf 0.49999999999999864688   \\\hline
\end{tabular}
}
\end{center}
\end{table}

\end{document}